\documentclass{amsart}
\usepackage{amssymb,amsmath,latexsym}

\input xy
\xyoption{all}

\font\tengoth=eufb10
\font\sevengoth=eufb7
\font\fivegoth=eufb5
\newfam\gothfam
\textfont\gothfam=\tengoth
\scriptfont\gothfam=\sevengoth
\scriptscriptfont\gothfam=\fivegoth
\def\goth{\fam\gothfam\tengoth}

\def\C{{\mathchoice {\setbox0=\hbox{$\displaystyle\mathrm C$}\hbox{\hbox
to0pt{\kern0.4\wd0\vrule height0.9\ht0\hss}\box0}}
{\setbox0=\hbox{$\textstyle\mathrm C$}\hbox{\hbox
to0pt{\kern0.4\wd0\vrule height0.9\ht0\hss}\box0}}
{\setbox0=\hbox{$\scriptstyle\mathrm C$}\hbox{\hbox
to0pt{\kern0.4\wd0\vrule height0.9\ht0\hss}\box0}}
{\setbox0=\hbox{$\scriptscriptstyle\mathrm C$}\hbox{\hbox
to0pt{\kern0.4\wd0\vrule height0.9\ht0\hss}\box0}}}}
\def\N{{\mathrm{I\!N}}} %natuerliche Zahlen

\def\Q{{\mathchoice {\setbox0=\hbox{$\displaystyle\mathrm
Q$}\hbox{\raise
0.15\ht0\hbox to0pt{\kern0.4\wd0\vrule height0.8\ht0\hss}\box0}}
{\setbox0=\hbox{$\textstyle\mathrm Q$}\hbox{\raise
0.15\ht0\hbox to0pt{\kern0.4\wd0\vrule height0.8\ht0\hss}\box0}}
{\setbox0=\hbox{$\scriptstyle\mathrm Q$}\hbox{\raise
0.15\ht0\hbox to0pt{\kern0.4\wd0\vrule height0.7\ht0\hss}\box0}}
{\setbox0=\hbox{$\scriptscriptstyle\mathrm Q$}\hbox{\raise
0.15\ht0\hbox to0pt{\kern0.4\wd0\vrule height0.7\ht0\hss}\box0}}}}
\def\Z{{\mathchoice {\hbox{$\sf\textstyle Z\kern-0.4em Z$}}
{\hbox{$\sf\textstyle Z\kern-0.4em Z$}}
{\hbox{$\sf\scriptstyle Z\kern-0.3em Z$}}
{\hbox{$\sf\scriptscriptstyle Z\kern-0.2em Z$}}}}

\def \Z {\mathbb Z}
\def \Hom {\rm Hom}
\def \FHom{\rm FHom}

\def \Ker {\rm Ker}
\def \Im {\rm Im}
\def \Ext {\rm Ext}
\newcounter{amoi}
\newtheorem{theorem}{Theorem}[subsection]
\newtheorem{lemma}[theorem]{Lemma}
\newtheorem{proposition}[theorem]{Proposition}%[section]
\newtheorem{corollary}[theorem]{Corollary}%[section]
\newtheorem{definition}[theorem]{Definition}%[section]

\def\hfl#1#2{\smash{\mathop{\hbox to 12 mm{\rightarrowfill}}
\limits^{\scriptstyle#1}_{\scriptstyle#2}}}

\begin{document}

\title{\large D\MakeLowercase{uality} \MakeLowercase{properties} \MakeLowercase{for} 
\MakeLowercase{quantum} \MakeLowercase{groups}}
\author{\large S\MakeLowercase{ophie} C\MakeLowercase{hemla}}

\maketitle
\centerline{UPMC Universit\'e Paris 6}
\centerline{UMR 7586}
\centerline{Institut de math\'ematiques}
\centerline{75005 Paris, France}
\centerline{schemla@math.jussieu.fr}

\vspace{2em}
{\bf Abstract :}
Some  duality properties for induced representations of
enveloping algebras involve the character $Trad_{\goth g}$. 
We extend them to deformation Hopf algebras
$A_{h}$ of a noetherian Hopf $k$-algebra $A_{0}$ satistying 
$Ext^{i}_{A_{0}}\left (k, A_{0}\right )=\{0\}$ except for $i=d$ where it is isomorphic to $k$. 
These
duality properties involve the  character of $A_{h}$ defined by right 
multiplication on the one dimensional free $k[[h]]$-module 
$Ext^{d}_{A_{h}}\left (k[[h]], A_{h} \right )$. In the case of quantized
  enveloping algebras, this character lifts the character 
$Trad_{\goth g}$. We also prove Poincar{\'e} duality for such
deformation Hopf algebras in the case where $A_{0}$ is of finite homological dimension. We explain the relation of our construction
with quantum duality.     

\section{Introduction}
 
In this article $k$ will be a field of characteristic $0$ and we set 
$K=k[[h]]$. 

Let $A_{0}$ be a noetherian algebra. 
We assume moreover that $k$ has a left
$A_{0}$-module structure such that there exists an integer $d$ 
satisfying 
$$\left \{ \begin{array}{l}
Ext_{A_{0}}^{i}\left ( k, A_{0} \right )=\{0\}\; {\rm if} \; i \neq d \\
Ext_{A_{0}}^{d} (k,A_{0}) \simeq k.  
\end{array}\right . $$ 

It follows from Poincar{\'e} duality that any finite dimensional Lie algebra
${\mathfrak g}$ verifies these assumptions. In this case $d=dim
{\mathfrak g}$ and the character defined by the right representation
of $U({\mathfrak g})$ on 
$Ext_{U({\mathfrak g})}^{dim {\mathfrak g}} (k,U({\mathfrak g}))$
is $Tr ad_{\mathfrak g}$ ([C1]). The algebra of regular 
fonctions on an affine algebraic Poisson  group and algebra of formal
power series also satisfy these hypothesis. 
Let $A_{h}$ be a deformation algebra of $A_{0}$. Assume that there exists
an $A_{h}$-module structure on $K$ that reduces modulo $h$ to the
$A_{0}$-module structure we started with.  
The following theorem
constructs a new character of  $A_{h}$, which will be denoted by 
$\theta_{A_{h}}$.\\

{\bf Theorem ~\ref{qtrad}}

{\it With the assumtions made above, one has :

a) $Ext^{i}_{A_{h}}(K, A_{h})=\{0\}$ is zero if $i\neq d$

b) $Ext^{d}_{A_{h}}(K, A_{h})$ is a free $K$-module of dimension
one. The right $A_{h}$-module structure given by right multiplication
lifts that of $A_{0}$ on $Ext_{A_{0}}^{d} (k,A_{0})$. }\\   

The right $A_{h}$-module $Ext^{d}_{A_{h}}(K, A_{h})$ will
be denoted by $\Omega_{A_{h}}$.  If there is an ambiguity, the integer $d$
will be written  $d_{A_{h}}$. 

Theorem $~\ref{qtrad}$  applies to universal quantum enveloping
algebras, quantization of affine algebraic
Poisson  groups and to quantum formal series Hopf algebras. 

Let ${\goth g}$ be a Lie bialgebra. Denote by   
$F[{\goth g}]$ the 
formal series Poisson  algebra $U({\goth g})^{*}$. 
If $U_{h}({\goth g}^{*})$ is a quantum enveloping algebra such that 
$U_{h}({\goth g}^{*})/hU_{h}({\goth g}^{*})$ is isomorphic to
$U({\goth g}^{*})$ as a coPoisson Hopf algebra, we show that one may construct a resolution
of  the trivial $U_{h}({\goth g}^{*})$-module $k[[h]]$ that lifts the Koszul
resolution of the trivial $U({\goth g}^{*})$-module $k$. 
If $F_{h}[{\goth g}]$ is a quantum formal
series algebras such that $F_{h}[{\goth g}]/hF_{h}[{\goth g}]$ is
isomorphic to $F[{\goth g}]$ as a Poisson Hopf algebra, we construct a
resolution of the trivial $F_{h}[{\goth g}]$-module that lifts the
Koszul resolution of the trivial $F[{\goth g}]$-module $k$ and that respects quantum duality
([Dr], [Ga]). This construction is not explicit but it allows to show that, if  
$F_{h}[{\goth g}]$ and $U_{h}({\goth g}^{*})$ are linked by quantum
duality, the following equality holds  $\theta_{F_{h}[{\goth g}]}=h\theta_{U_{h}({\goth g}^{*})}$.\\

As an application of theorem ~\ref{qtrad}, we show Poincar{\'e}
duality : \\

{\bf Theorem ~\ref{Poincare}}

{\it We make the same assumtions as above. 
Let $M$ be an $A_{h}$-module. 
Assume that $K$ is an $A_{h}$-module of finite projective dimension. One has an isomorphism
of $K$-modules for all integer $i$ :}\\

\centerline{$Ext^{i}_{A_{h}}\left ( K, M \right ) \simeq 
Tor^{A_{h}}_{d_{A_{h}}-i}\left ( \Omega_{A_{h}}, M \right ).$}
\vspace{2em}
{\it From now on, we assume that $A_{h}$ is a deformation Hopf algebra.}\\

Brown and Levasseur ([B-L]) and Kempf ([Ke]) had shown that, in the
semi-simple context, the Ext-dual of a Verma module is a Verma
module. In [C1], we have extended this result to the Ext-dual of an
induced representation of any Lie superalgebra. In this article, we show
that this result can be generalized to quantum groups provided that
the quantization is functorial. Such a functorial quantization has
been constructed by Etingof and Kazdhan ([E-K1], [E-K2], [E-K3], [E-S]). 
As the result holds for quantized universal enveloping algebras, for 
quantized functions algebras and for quantum formal series Hopf
algebras, we state it in the more general setting of Hopf
algebras. \\

{\bf Corollary ~\ref{corollaire Ext-dual}}

{\it Let $A_{h}$ (respectively $B_{h}$)
be a topological Hopf deformation of $A_{0}$ (respectively $B_{0}$). 
We assume that there exists a
morphism of Hopf algebras from $B_{h}$ to $A_{h}$ such that $A_{h}$ is
a flat $B_{h}^{op}$-module. We also assume that $B_{h}$ satisfies the
condition of the theorem ~\ref{qtrad}. Let $V$ be a $B_{h}$-module which is
a free finite dimensional $K$-module. Then 

a) $Ext^{i}_{A_{h}} \left ( 
A_{h}{\displaystyle \mathop \otimes_{B_{h}}}V, A_{h} \right )$ is
 $\{0\}$ if $i$ is different from $d_{B_{h}}$.

b) The right $A_{h}$-module 
 $Ext^{d_{B_{h}}}_{A_{h}} \left ( 
A_{h}{\displaystyle \mathop \otimes_{B_{h}}}V, A_{h} \right )$ is
isomorphic to 
$\left (  \Omega_{B_{h}}\otimes V^{*}\right ){\displaystyle \mathop \otimes_{B_{h}}}A_{h}$} where 
$ \Omega_{B_{h}}\otimes V^{*}$ is endowed with the following
right  $B_{h}$-module structure :
$$\begin{array}{l}
\forall u \in B_{h}\, \forall f \in V^{*}, \;
\forall \omega \in \Omega_{B_{h}},\\
(\omega \otimes f)\cdot u =
 \lim\limits_{n\to +\infty}
{\displaystyle \mathop \sum_{j}}
\theta_{B_{h}} (u'_{j,n})\omega  \otimes 
f \cdot S_{h}^{2}(u''_{j,n}) \\
\Delta (u)=\lim\limits_{n\to +\infty}
{\displaystyle \mathop \sum_{j}u'_{j,n}\otimes u''_{j,n}}.
\end{array}$$
$S_{h}$ being the antipode of $B_{h}$. \\

{\bf Proposition ~\ref{duality property}}
{\it Let  $A_{h}$ be a Hopf deformation of $A_{0}$,  $B_{h}$ be a
Hopf  deformation of $B_{0}$ and  
$C_{h}$ be a Hopf  deformation of $C_{0}$. We assume that there
exists a morphism of Hopf algebras from $B_{h}$ to $A_{h}$ 
and a morphism of Hopf algebras from $C_{h}$ to $A_{h}$ 
such that $A_{h}$
is a flat  $B_{h}^{op}$-module and a flat $C_{h}^{op}$-module. We also
assume that $B_{h}$ and $C_{h}$ satisfies the hypothesis of 
theorem ~\ref{qtrad}. 
Let $V$ (respectively $W$) be a 
$B_{h}$-module (respectively $C_{h}$-module) which is  a free finite 
dimensional K-module. Then, for all integer $n$, one has an
isomorphism 

$$\begin{array}{l}
Ext^{n +d_{B_{h}}}_{A_{h}}\left ( 
A_{h}{\displaystyle \mathop \otimes_{B_{h}}} V, 
A_{h}{\displaystyle \mathop \otimes_{C_{h}}} W
\right )\\
\simeq 
Ext^{n +d_{C_{h}}}_{A_{h}}\left ( 
\left ( \Omega_{C_{h}}\otimes W^{*}\right )
{\displaystyle \mathop \otimes_{C_{h}}} A_{h}, 
\left (  \Omega_{B_{h}}\otimes V^{*}\right )
{\displaystyle \mathop \otimes_{B_{h}}}A_{h}  \right ).
\end{array} $$}
The right $B_{h}$ (respectively $C_{h}$)-module structure on 
$ \Omega_{B_{h}}\otimes V^{*}$ (respectively $ \Omega_{C_{h}}\otimes W^{*}$ ) are as in 
 Corollary ~\ref{corollaire Ext-dual}.\\

{\bf Remark : } 

Proposition ~\ref{duality property} is already known in
the case where $\mathfrak g$ is a Lie algebra, ${\mathfrak h}$ and 
${\mathfrak k}$ are Lie subalgebras of ${\mathfrak g}$, 
$A_{h}$, $B_{h}$ and  
 $C_{h}$ are their corresponding enveloping algebras. 
In this case one has $d_{B_{h}}=dim{\mathfrak h}$ and 
$d_{C_{h}}=dim{\mathfrak k}$. More
 precisely : 

 Generalizing a result of G. Zuckerman ([B-C]), A. Gyoja ([G]) proved a
 part of this theorem (namely the case where 
${\mathfrak h}={\mathfrak g}$ and $n=dim{\mathfrak h}=dim {\mathfrak
  k}$) under the assumptions that ${\mathfrak g}$ is split semi-simple
and ${\mathfrak h}$ is a parabolic subalgebra of ${\mathfrak g}$. 
D.H Collingwood and B. Shelton ([C-S]) also proved a duality of this
type (still under the semi-simple hypothesis) 
but in a slighly different context.  

M. Duflo [Du2] proved proposition ~\ref{duality property} for a
${\mathfrak g}$ general Lie algebra, ${\mathfrak h}={\mathfrak k}$,
$V=W^{*}$ being one dimensional representations.

Proposition ~\ref{duality property} is proved in full generality in the
context of Lie superalgebras in [C1].\\

Wet set $A_{h}^{e}=A_{h}\otimes A_{h}^{op}$. Using the properties as a Hopf algebra (as in [C2]), we show that all the 
$Ext^{i}_{\widehat{A_{h}^{e}}}(A_{h}, A_{h}\widehat{\displaystyle \mathop \otimes _{k[[h]]}}A_{h})$'s 
are zero except one. More precisely :

{\bf Proposition \ref{Hochschild}}
{\it Assume that $A_{h}$ satisfies the conditions of the theorem \ref{qtrad}. Assume moreover that 
$A_{0}\otimes A_{0}^{op}$ is noetherian. Consider 
$A_{h}\widehat{\displaystyle \mathop \otimes _{k[[h]]}}A_{h}$ with the following 
$\widehat{A_{h}^{e}}$-module structure :
$$\forall (\alpha , \beta , x ,y ) \in A_{h}, \;\;\;\alpha \cdot (x \otimes y) \cdot \beta = 
\alpha x \otimes y \beta .$$

a) $HH^{i}_{A_{h}}(A_{h}
\widehat{\displaystyle \mathop \otimes _{k[[h]]}}A_{h})$ is zero if $i\neq d_{A_{h}}$.

b)  The $\widehat{A_{h}^{e}}$-module  $HH^{d_{A_{h}}}_{A_{h}}
(A_{h}\widehat{\displaystyle \mathop \otimes _{k[[h]]}}A_{h})$ is isomorphic to 
$\Omega_{A_{h}} \otimes A_{h}$ with the following $\widehat{A_{h}^{e}}$-module structure : 
$$\forall (\alpha , \beta , x ) \in A_{h}, \;\;\;\ \\
\alpha \cdot (\omega \otimes x) \cdot \beta = 
\omega \theta_{A_{h}}(\beta '_{i})  \otimes S(\beta_{i}'')
 x S^{-1}( \alpha) $$
 where $\alpha={\displaystyle \sum_{i}}\alpha_{i}'\otimes \alpha''_{i}$ (to be taken in the topological sense)}

From this, as in [VdB], we deduce a relation between Hochschild homology and Hochschild cohomology
for the ring $A_{h}$.\\

{\bf Acknowledgments :}

I am grateful to B. Keller, D. Calaque, B. Enriquez and V. Toledano for 
helpful discussions.
 
\section{Graded linear algebra}

In this section, we fix notation about graded linear algebra. A graded
$k$-algebra $GA$ is the data of a $k$-algebra with unit and a family of $k$-vector spaces  
$(G_{t}A)_{t \in \Z}$ of $A$ satisfying~:
$$\begin{array}{l}
a. \;\;\; A={\displaystyle \oplus _{t \in \Z}}G_{t}A \\

b.\;\;\; 1 \in G_{0}A \\

c.\;\;\; G_{t}A \cdot G_{l}A \subset G_{t+l}A.
\end{array}$$
We will also assume that $G_{t}A=0$ for $t<0$.

A graded $GA$-module $GM$ is the data of a $GA$-module and a family of
$k$-vector space $\left ( G_{t}M \right )_{t \in \Z}$ of $GM$ 
such that 
$$\begin{array}{l}
GM={\displaystyle \oplus_{t \in \Z}}G_{t}M \\
G_{t}A \cdot G_{l}M \subset G_{t+l}M   
\end{array}$$
We will always also assume that $G_{t}M=0$ if $t<<0$. 

Let $GM$ and $GN$ be two graded $GA$-modules. A morphism of graded
$GA$-modules from $GM$ to $GN$ is a morphism of $GA$-modules 
$f : GM \to GN$ such that $f(G_{t}M)\subset G_{t}N$. The group of
morphisms of graded $GA$-modules from $GM$ to $GN$ will be denoted 
$\Hom_{GA}(GM ,GN)$. With this notion of morphisms,  the category of
graded $GA$-modules is abelian. Thus it is suitable for homological
algebra. 

For $r \in {\mathbb Z}$ and any graded $GA$-module $GM$, we define the
shifted graded $GA$-module $GM(r)$ to be the $GA$-module $GM$ endowed with
the grading defined by 
$$\forall t \in \Z , \;\;\; G_{t}M(r)=G_{t+r}M.$$
Let us denote $\underline{\rm Hom}_{GA}(GM, GN)$ the graded group defined by
setting 
$$G_{t}\underline{\Hom }_{GA}(GM ,GN)=Hom_{GA}(GM, GN (t)).$$
The i$^{th}$ right derived functor of the functor $\underline{\rm Hom}_{GA}(-,N)$ will be
denoted $\underline{\rm Ext}^{i}_{GA}(-,N)$. 

A graded $GA$-module $GL$ is finite free if there are integers 
$d_{1},d_{2},\dots ,d_{n}$ such that 
$$GL\simeq {\displaystyle \mathop \oplus_{i=1}^{n}}GA(-d_{i}).$$
A graded $GA$-module $GM$ is of finite type if there exists a finite
free graded $GA$-module $GL$ and an exact sequence in the category of
graded $GA$-modules  
$$GL \to GM \to 0.$$ 
This means that there are homogeneous elements $m_{1}\in G_{d_{1}}M,
\dots ,m_{n} \in G_{d_{n}}M$ such that any $m \in G_{d}M$ may be
written as 
$$m={\displaystyle \sum_{i=1}^{n}}a_{d-d_{i}}m_{d_{i}}$$
where $a_{d-d_{i}} \in G_{d-d_{i}}A.$

A graded ring $GA$ is noetherian if any graded $GA$-submodule of a
graded $GA$-module of finite type is of finite type. \\

In the sequel, all the $GA$-modules we will consider will be graded so
that we will say "$GA$-module" for "graded $GA$-module".

\section{Decreasing filtrations}

In this section, we give results about decreasing filtrations. These 
results are proved in [Schn] in the framework of increasing
filtrations. 
For the sake of completeness, we give detailled proofs of
the results even if most of our proofs are obtained by adjusting those
of Schneiders.  \\

We will consider a $k$-algebra endowed with a decreasing filtration 
$...F_{t+1}A\subset F_{t}A \subset \dots \subset F_{1}A\subset
F_{0}A=A$. The order of an element $a$, $o(a)$,  is the biggest $t$ such that 
$a \in F_{t}A$. The principal symbol of $a$ is the image of $a$ in 
$F_{o(a)}/F_{o(a)+1}$. It will be denoted by $[a]$.  

A filtered module over $FA$ is the data of an $A$-module $M$ and a
family $\left ( F_{t}M \right )_{t \in \Z }$ of $k$-subspaces such
that 

\begin{itemize}
\item ${\displaystyle \bigcup _{t \in \Z}}F_{t}M =M$\\
\item $F_{t+1}M \subset F_{t}M$
\item $F_{t}A \cdot F_{l}M \subset F_{t+l}M$
\end{itemize} 

We will assume that $F_{t}M=M$ for $t<<0$. 
We have the notion of principal symbol. 
We endow such a module with the topology for which a basis a
neighborhoods is $\left ( F_{t}M \right )_{t \in \Z}$. The topological
space $M$ is Hausdorff if and only if 
${\displaystyle \cap_{t\in \Z}}F_{t}M=\{0\}$. If $M$ is Hausdorff, the
topology defined by the filtration is defined by the following metric 
$$\begin{array}{l}
\forall (x,y) \in FM ,\; d(x,y)=\mid \mid x-y \mid \mid \;\;
{\rm with}\;\; \\ 
\mid \mid x-y \mid \mid = 2^{-t} \;\;{\rm where}
\; t=Sup \{j \in \Z \mid x-y \in F_{j}M \}
\end{array}$$
Note that $M$ is Hausdorff if and only if the natural map from $M$ to 
${\displaystyle \varprojlim_{t \in \Z }}{\displaystyle \frac
  {M}{F_{t}M}}$ is injective. The metric
space $(M,d)$ is  complete   if and only if the natural map from $M$ to
${\displaystyle \varprojlim_{t \in \Z }} 
{\displaystyle \frac {M}{F_{t}M}}$ is an isomorphism.\\

{\bf Example : }

Let $k$ be a field and set $K=k[[h]]$. 
If $V$ is a $K$-module, it is endowed with the following decreasing
filtration $\dots \subset h^{n}V \subset h^{n-1}V \subset \dots
\subset hV\subset V$.\\ 
The topology induced by this filtration is the $h$-adic topology.\\

Recall the following result :

\begin{lemma}
Let $N$ be a Hausdorff filtered module. Let $P$ be a submodule of $N$ which is
closed in $N$. Let $p$ the canonical projection from $N$ to $N/P$. 

a) The topology defined by the filtration $p(F_{t}N)$ on $N/P$ is the
quotient topology. $N/P$ is Hausdorff and its  topology is defined by the 
distance $d(\bar{x},\bar{y})=\mid\mid \bar{x}-\bar{y}\mid \mid $ where 
$$\mid \mid \bar{x}\mid \mid = 
Inf\{\mid \mid a \mid \mid , a \in \bar{x} \}$$

b) If $N$ is complete, then $N/P$ is complete for the quotient
topology. 

\end{lemma}

{\it Proof of the lemma :}

a) As $P$ is closed in $N$, then $\bar{0}$ is closed in $N/P$. Thus, its
complement in $N/P$, $U$, is open. Let $\bar{x}$ an element of $N/P$
different from $\bar{0}$. As $U$ is open, there exists $n \in \N$
such that $\bar{x}\in p(F_{n}N)+\bar{x} \subset U$. Hence 
$\bar{x} \notin p(F_{n}N)$ and we have proved that 
${\displaystyle \bigcap _{n \in \N}}p(F_{n}N)=\{\bar{0}\}$. Hence
$N/P$ is Hausdorff. 
It is easy to check that the open ball of center $0$ and radius $2^{-t}$
in $N/P$ for the distance  defined is $p(F_{t+1}N)$.

b)  we refer to [Schw] p 245.   
$\Box$. \\

Let $FM$ and $FN$ be two filtered $FA$-modules. 
A filtered morphism $Fu : FM \to FN$ is a morphism $u :M \to N$ 
of the underlying $A$-modules such that $u(F_{t}M) \subset F_{t}N$.
It is continuous if we endow $M$ and $N$ with the topology defined by
the filtrations.   
Denote by $F_{t}u$ the morphism $u_{\mid F_{t}M} : F_{t}M \to
F_{t}N$. Denote by 
$Hom_{FA}(FM, FN)$ the group of filtered morphisms from $FM$ to $FN$.  
The kernel of $Fu$ is the kernel of $u$ filtered by the family 
$KerFu \cap F_{t}M$. If $M$ is complete and $N$ is Hausdorff, then  
 $Ker Fu$,  endowed with the induced
topology is  complete. \\

To a filtered ring $FA$ is associated a graded ring $GA$ defined by 
$$GA={\displaystyle \mathop \oplus _{t \in \N}}G_{t}A \;\;{\rm with
}\;\;G_{t}A=F_{t}A /F_{t+1}A$$
the multiplication being induced by that of $FA$. 
To a filtered $FA$-module $FM$ is  associated  a graded $GA$-module 
$GM$ defined by setting  
$$GM={\displaystyle \mathop \oplus _{t \in \Z}}G_{t}M \;\;{\rm with
}\;\;G_{t}M=F_{t}M /F_{t+1}M$$
the action of $GA$  on $GM$ being induced by that of $FA$. 
If $x$ is in $F_{t}M$, we will write $\sigma_{t}(x)$ for the class of
$x$ in $F_{t}M/F_{t+1}M$. 
A filtered
morphism of $FA$-modules $Fu : FM \to FN$ induces a morphism 
of abelian groups $G_{t}u : G_{t}M \to G_{t}N$ and a  morphism of 
$GA$-modules $Gu : GM \to GN$.\\

An arrow $Fu : FM \to FN$ is strict if it satisfies 
$u(F_{t}M)=u(M)\cap F_{t}N$. \\

An exact sequence of $FA$-modules is a sequence 
$$FM \buildrel{Fu}\over \longrightarrow FN  
\buildrel{Fv} \over \longrightarrow FP $$ 
such that $Ker F_{t}v=ImF_{t}u$.  It follows from this definition that 
$Fu$ is strict. if moreover $Fv$ is strict, we say that it is a strict
exact sequence. \\

\begin{proposition}~\label{strict exact sequences} 
a) Consider $Fu : FM \to FN$ and $Fv : FN \to FP$ two filtered
$FA$-morphisms such that $Fv \circ Fu=0$. If the sequence 
$$FM \buildrel{Fu}\over \longrightarrow FN  
\buildrel{Fv} \over \longrightarrow FP $$ 
is strict exact, then 
$$GM \buildrel{Gu}\over \longrightarrow GN  
\buildrel{Gv} \over \longrightarrow GP $$ 
is exact. 

b) Conversely, assume that $FM$ is  
complete for the topology 
defined by the filtration and that $FN$ is Hausdorff for the topology
defined by the filtration. If the sequence 
$$GM \buildrel{Gu}\over \longrightarrow GN  
\buildrel{Gv} \over \longrightarrow GP $$ 
is exact, then the sequence 
$$FM \buildrel{Fu}\over \longrightarrow FN  
\buildrel{Fv} \over \longrightarrow FP $$ 
is strict exact.
\end{proposition}

{\it Proof of the proposition :}

a) Let $n_{t} \in G_{t}N$ be such that $G_{t}v(n_{t})=0$. There is
$n'_{t} \in F_{t}N$ such that $n_{t}=\sigma_{t}(n'_{t})$. Hence 
$v(n'_{t})\in F_{t+1}P$. Since $Fv$ is strict, we find $n''_{t+1}\in F_{t+1}N$ 
such that $v(n''_{t+1})=v(n'_{t})$. Then $v(n'_{t}-n''_{t+1})=0$
and there is $m_{t}\in F_{t}M$ such that
$u(m_{t})=n'_{t}-n''_{t+1}$. This shows that 
$$G_{t}u\left ( \sigma_{t}(m_{t}) \right )= \sigma_{t}(n'_{t})=n_{t}.$$

b ) Let us prove that $Fv$ is strict. Assume that $p_{t} \in F_{t}P \cap
Imv$. Let $l$ be the biggest integer such that $p_{t}=v(n_{l})$ with 
$n_{l} \in F_{l}N$. We need to show that $l \geq t$. Assume that $l
< t$. One has 
$$G_{l}v\left ( \sigma_{l}(n_{l})\right )=
\sigma_{l}\left ( v(n_{l})\right )=\sigma_{l}(p_{t})=0.$$
Hence $\exists m_{l} \in F_{l}M$ such that 
$G_{l}u\left ( \sigma_{l}(m_{l}) \right )=\sigma_{l}(n_{l})$. Thus we
have 
$$n_{l}-u(m_{l}) \in F_{l+1}M\;\; {\rm and} \;\; 
v\left ( n_{l}-u(m_{l}) \right ) =p_{t}$$ 
which contredicts the definition of $l$. \\

Let us prove that $KerF_{t}v=ImF_{t}u$. Let $n_{t} \in KerF_{t}v$.
One has :  
$G_{t}(v)\left ( \sigma_{t}(n_{t}) \right )=0$. Hence there exists
$m_{t}$ in $F_{t}M$ such that 
$$\sigma _{t}(n_{t})=G_{t}(u)\left ( \sigma_{t}(m_{t})\right ).$$
Hence $n_{t}-u(m_{t}) \in Ker Fv \cap F_{t+1}N$. 
We can reproduce the previous reasoning to $n_{t}-u(m_{t})$ and
produce an element $m_{t+1}$ in $F_{t+1}M$ such that 
$n_{t}-u(m_{t}+m_{t+1}) \in Ker Fv \cap F_{t+2}N$. The sequence 
${\mathcal U}_{p}={\displaystyle \sum_{l=0}^{p}m_{t+l}}$ is a Cauchy
sequence, hence it converges and 
$n_{t}=u \left ( {\displaystyle \sum_{l=0}^{\infty}m_{t+l}}\right )$.
$\Box$\\

\begin{corollary}
Let $FA$ be a filtered $k$-algebra and let $FM$ and $FN$ two
$FA$-modules. Let   
$Fu : FM \to FN$  be a morphism of
$FA$-modules. Then $G \Ker Fu \subset \Ker GFu$ and 
$\Im G Fu \subset G\Im Fu$. 
Assume moreover that  $FM$ is complete and $FN$ is Hausdorff, 
then the following conditions are equivalent : 

(a) $Fu$ is strict

(b) $G \Ker Fu = \Ker GFu$

(c) $\Im G Fu = G\Im Fu$. 
\end{corollary}

{\it Proof :}

 One has : 

$$\begin{array}{l}
F_{t}\Ker u=\Ker u\cap F_{t}M\\ 
F_{t}\Im u=\Im u\cap F_{t}N\\
G_{t}\Ker u={\displaystyle \frac{F_{t}M \cap \Ker u}
{F_{t+1}M \cap \Ker u}}\\
\Ker G_{t} u={\displaystyle \frac {F_{t}M\cap u^{-1}(F_{t+1}N)}
{F_{t+1}M\cap u^{-1}(F_{t+1}N)}}\\
\Im G_{t}u={\displaystyle \frac{u(F_{t}M)}
{F_{t+1}N\cap u(F_{t}M)}}\\
G_{t}\Im u={\displaystyle \frac{\Im u\cap F_{t}N}
{\Im u\cap F_{t+1}M}}
\end{array}$$
The second part of the corollary follows from applying the previous 
proposition to the
strict exact sequence $FM \to \Im u \to 0$.

Indeed 
$Fu$ is strict if and only the following sequence  
$$FM \buildrel {Fu} \over \longrightarrow  \Im u \to 0$$
is a strict exact sequence  of $FA$-modules  when $\Im u$ 
is endowed with the induced topology.   Then we apply ~\ref{strict exact sequences} . \\

Let us recall this well known result about complexes of filtered
modules. \\

\begin{proposition}~\label{spectral sequence}
Let $(M^{\bullet}, d^{\bullet})$ be a complex of complete 
$FA$-modules. $H^{i}(M^{\bullet})$ is filtered as follows 
$F_{t}H^{i}(M^{\bullet})=
{\displaystyle \frac{Kerd_{i}\cap F_{t}M^{i}+Imd_{i-1}}{Imd_{i-1}}}
\simeq 
{\displaystyle \frac{Kerd_{i}\cap F_{t}M^{i}}{Imd_{i-1}\cap
    F_{t}M^{i-1}}}$.  If  $d_{i}$ and $d_{i-1}$ are strict, then 
    $GH^{i}(M^{\bullet})$ is isomorphic to  $H^{i}(G M^{\bullet})$
 %The $GA$-module $GH^{i}(M^{\bullet})$ is a subfactor of  $H^{i}(G M^{\bullet})$. 
\end{proposition}

{\it Proof of the proposition ~\ref{spectral sequence}: }
 
We consider the following exact sequence 
$$0\to Imd_{i-1}\to Kerd_{i} \buildrel{p}\over \longrightarrow
H^{i}M^{\bullet}\to 0.$$
we endow  $Kerd_{i}$ and $Imd_{i-1}$ with the induced filtration.   
One has 
$$\begin{array}{l}
F_{t}Kerd_{i}=Kerd_{i}\cap F_{t}M^{i}\\ 
F_{t}Imd_{i-1}=Imd_{i-1}\cap F_{t}M^{i}\\
p(F_{t}Kerd^{i})=
{\displaystyle
  \frac{Kerd_{i}\cap
    F_{t}M^{i}+Imd_{i}}{Imd_{i}}}=F_{t}H^{i}(M^{\bullet}).
\end{array}$$
The exact sequence above is strict exact. It stays exact if one takes
the graded modules. Thus, we have the following exact sequence of 
$GA$-modules
 $$0\to GImd_{i-1}\to GKerd_{i} \buildrel{p}\over \longrightarrow
G H^{i}M^{\bullet}\to 0.$$ 
%One has 
%$$\begin{array}{l}
%G_{k}Kerd_{i}={\displaystyle \frac{F_{k}M^{i} \cap Kerd_{i}}
%{F_{k+1}M^{i} \cap Kerd_{i}}}\\
%Ker G_{k} d_{i}=
%{\displaystyle \frac {F_{k}M^{i}\cap d^{-1}_{i}(F_{k+1}M^{i+1})}
%{F_{k+1}M^{i}\cap d_{i}^{-1}(F_{k+1}M^{i+1})}}\\
%ImG_{k}d_{i-1}={\displaystyle \frac{d_{i-1}(F_{k}M^{i-1})}
%{F_{k+1}M^{i}\cap d_{i-1}(F_{k}M^{i-1})}}\\
%G_{k}Imd_{i-1}={\displaystyle \frac{Imd_{i-1}\cap F_{k}M^{i}}
%{Imd_{i-1}\cap F_{k+1}M^{i}}}
%\end{array}$$ 
%We have the following sequence of inclusions 
%$$ImGd_{i-1}\subset GImd_{i-1} \subset GKerd_{i} \subset
%KerGd_{i}.$$
%Hence
%$$H={\displaystyle \frac {GImd_{i-1}}{ImGd_{i-1}}} 
%\subset G={\displaystyle \frac{GKerd_{i}}{ImGd_{i-1}}} 
%\subset {\displaystyle \frac{KerGd_{i}}{ImGd_{i-1}}}
%\simeq  H^{i}(GM^{\bullet}).$$
%and $GH^{i}(M^{\bullet})\simeq 
%{\displaystyle \frac{G}{H}}$. 
Then $G H^{i}(M^{\bullet})\simeq {\displaystyle \frac{GKerd_{i} }{GImd_{i-1} }} \simeq 
{\displaystyle \frac{KerGd_{i} }{ImGd_{i-1} }}\simeq H^{i}(GM^{\bullet})$. 
This finishes the proof of the proposition. 
$\Box$\\

{\bf Remark :}

The isomorphism from $G_{t} H^{i}(M^{\bullet})$ to $ H^{i}(G_{t}M^{\bullet})$ is given by 
$$\begin{array}{rcl}
G_{t} H^{i}(M^{\bullet} )&\to & H^{i}(G_{t}M^{\bullet})\\
\sigma_{t}cl(x)& \mapsto & cl(\sigma _{t}(x)).\\
\end{array}$$\\

For any $r \in \Z$ and for any $FA$-module $FM$, we define the
shifted module $FM(r)$ as the module $M$ endowed with the filtration 
$\left ( F_{t+r}M \right )_{t \in \Z}$.\

An $FA$-module module is  finite free if it is isomorphic to an 
$FA$-module of the type 
${\displaystyle \oplus_{i=1}^{p}FA(-d_{i})}$ where 
$d_{1},\dots , d_{p}$ are integers. An $FA$-module $FM$ is of finite
type if there exists a strict epimorphism $FL \to FM$ where  $FL$ is a
finite free $FA$-module. This means that we can find 
$m_{1} \in F_{d_{1}}M, \dots, m_{p}\in F_{d_{p}}M$ such that any 
$m \in F_{d}M$ may be written as 
$$m={\displaystyle \sum_{i=1}^{p}a_{d-d_{i}}}m_{i}$$   
where $a_{d-d_{i}} \in F_{d-d_{i}}A$. 

\begin{proposition}
Let $FA$ be a filtered $k$-algebra and $FM$ be an $FA$-module. 

a) If $FM$ is an $FA$-module of finite type generated by 
$(s_{1},\dots , s_{r})$ then  $GM$ is a
$GA$-module of finite type generated by $([s_{1}],\dots , [s_{r}])$. 
Conversely, assume that $FA$ is  complete for the topology
given by the filtration and that  $FM$ is a $FA$-module 
which is Hausdorff  for the
topology defined by the filtration. If $GM$ is a 
$GA$-module of finite type generated by $([s_{1}],\dots , [s_{r}])$,
then $FM$ is an $FA$-module of finite type generated by 
$(s_{1},\dots , s_{r})$

b) If $FM$ is a finite free $FA$-module, then  $GM$ is a finite
free $GA$-module. Conversely, assume that $FA$ is  complete for the topology
given by the filtration and $FM$ is a $FA$-module Hausdorff for the
topology defined by the filtration. If $GM$ is a finite free
$GA$-module, then $FM$ is a finite free $FA$-module.
\end{proposition}

{\it Proof of the proposition :}

a) If $FM$ is an $FA$-module of finite type, then there is a strict  exact
sequence ${\displaystyle \oplus_{i=1}^{N}FA(-d_{i}) \to FM \to 0}$. 
If we apply
proposition 3.0.2, we see  that $GM$ is a $GA$-module of finite type. 
Conversely, assume that $GM$ is a $GA$-module of finite type 
generated by $\sigma_{1}=[s_{1}], \dots , \sigma_{r}=[s_{r}]$. Assume
that $s_{i} \in F_{d_{i}}M - F_{d_{i}-1}M$. Let $x$ in
$F_{n}M$. There exists $a_{i,0} \in G_{n-d_{i}}A$ such that 
$$\sigma_{n}(x)={\displaystyle \mathop
  \sum_{i=1}^{r}}a_{i,0}\sigma_{i}.$$
Let $\alpha_{i,0} \in F_{n-d_{i}}A$ such that 
$\sigma_{n-d_{i}}(\alpha_{i,0})=a_{i,0}$. We have 
$$x - {\displaystyle \mathop
  \sum_{i=1}^{r}}\alpha_{i,0}s_{i} \in F_{n+1}M.$$
Reasoning in the same way, one can construct 
$\alpha_{i,1}\in F_{n-d_{i}+1}A$ such that 
 $$x - {\displaystyle \mathop
  \sum_{i=1}^{r}}(\alpha_{i,0}+\alpha_{i,1})s_{i} \in F_{n+2}M $$
Going on that way, we construct an element 
${\displaystyle \mathop
  \sum_{j=1}^{\infty}}\alpha_{i,j}$ in $F_{n-d_{i}}A$
such that 
$$x =
{\displaystyle \mathop  \sum_{i=1}^{r}}
\left ( \sum_{j=1}^{\infty}\alpha_{i,j} \right ) s_{i}.$$ Hence $FM$
is a finite type $FA$-module. \\

b) apply proposition 3.0.2.$\Box$\\

\begin{definition}
A filtered $k$-algebra is said to be (filtered) noetherian if it satisfies
one of the following equivalent conditions :
\begin{itemize}
\item Any filtered submodule (not necessarily a strict submodule) of a
  finite type $FA$-module is of finite type
\item Any filtered ideal (not necessarily a strict ideal) of $FA$ is
  of finite type.  
\end{itemize}
\end{definition}

\begin{proposition}
Let $FA$ be a filtered complete $k$-algebra and denote by $GA$ its associated graded
algebra. If $GA$ is graded noetherian, then $FA$ is filtered noetherian.  
\end{proposition}

{\it Proof of the proposition :} 

We assume that GA is a noetherian algebra. 
We need to prove that a filtered submodule $FM'$ of a finitely
generated 
$FA$-module $FM$ is finitely generated. \\

{\it First we assume that $FM$ is Hausdorff}. For this case, we
reproduce the proof of [Sch]. 

If $FM'$ is strict, then the associated $GA$-module $GM'$ is a
submodule of the $GA$-module $GM$ associated to $FM$. Since $GA$ is
noetherian and $GM$ is  finitely generated  so is $GM'$ and the
conclusion follows. 

To prove the general case, we may assume that the image of the
inclusion $FM' \to FM$ is equal to $FM$. In this case, using a finite
systeme of generators of $FM$, it is easy to find an integer $l$ such
that 
$$F_{t}M' \subset F_{t}M \subset F_{t-l}M'.$$ 
We will prove the result by induction on $l$.

For $l=1$, let us introduce the auxiliary $GA$-modules 
$$\begin{array}{l}
GK_{0}={\displaystyle \oplus _{t \in \Z}}F_{t}M'/F_{t+1}M\\
GK_{1}={\displaystyle \oplus _{t \in \Z}}F_{t}M/F_{t}M'
\end{array}$$
These modules satisfy the exact sequences 
$$0 \to GK_{0} \to GM \to GK_{1}\to 0 $$
$$0 \to GK_{1}(1)\to GM'\to GK_{0}\to 0 $$
Since $GM$ is a finite type $GA$-module, so are $GK_{0}$ and
$GK_{1}$. Hence $GM'$ is also finitely generated and the conclusion
follows. 

For $l>1$, we define the auxiliary $FA$-module $FM''$ by setting 
$$F_{t}M''=F_{t+1}M+F_{t}M'.$$
Since we have 
$$F_{t}M'' \subset F_{t}M\subset F_{t-1}M''$$
the preceeding discussion shows that $FM''$ is finitely
generated. Moreover 
$$F_{t}M' \subset F_{t}M''\subset F_{t-(l-1)}M'$$
and the conclusion follows from the induction hypothesis. \\

{\it We no longer assume that $FM$ is Hausdorff}

As $FM$ is a finite type $FA$-module, there exists a strict exact
sequence 
$$FL={\displaystyle \mathop \oplus _{i=1}^{n}FA (-d_{i})} 
\buildrel {p}\over \longrightarrow FM \to 0.$$
We will denote by $p_{t}$ the map from $F_{t}L$ to $F_{t}M$ induced by
$p$. As $p$ is strict, the map $p_{t}$ is surjective. 
Let $FM'$ be a submodule (not necessarily strict) of $FM$. Then 
$p^{-1}(FM')$ is an $FA$-submodule of $FL$ if we endow it with the
filtration 
$$F_{t}\left [ p^{-1}(M')\right ]=
p_{t}^{-1}(F_{t}M')=p^{-1}(F_{t}M')\cap F_{t}L .$$
As $FL$ is Hausdorff, we know from the first part of the proof that
the 
$FA$-module $p^{-1}M'$ is finite type. Hence there exist $\alpha_{1}
\in F_{\delta_{1}}\left [ p^{-1}M'\right ], \dots , \alpha_{p}  \in
F_{\delta_{p}}\left [ p^{-1}M'\right ]$ 
such that any $x$ of $F_{d}\left [ p^{-1}M'\right ] $ can be written
$$x={\displaystyle \mathop
  \sum_{i=1}^{p}}a_{d-\delta_{i}}\alpha_{i}\;\;
{\rm with}\;\;a_{d-\delta_{i}}\in F_{d-\delta_{i}}A.$$
Let $y$ in $F_{d}M'$. As $p$ is strict, there exist 
$x \in F_{d}\left [ p^{-1}M'\right ]$ such that $y=p(x)$. 
Then $y$ can be written 
$$y={\displaystyle \mathop
  \sum_{i=1}^{p}}a_{d-\delta_{i}}p( \alpha_{i} )\;\;
{\rm with}\;\;a_{d-\delta_{i}}\in F_{d-\delta_{i}}A.$$
We have proved that $FM'$ is a finite type $FA$-module $\Box$.\\

\begin{proposition}~\label{noetherian}
Assume that $FA$ is noetherian  for the topology
given by the filtration.
Any $FA$-module of finite type has an infinite resolution by finite
free $FA$-modules i.e  there is an exact sequence
$$\dots \to FL_{s}\to FL_{s-1}\to \dots \to FL_{0}\to FM\to 0$$
where each $FL_{s}$ is a finite free $FA$-module. 
\end{proposition}

{\bf Remark :}

For such a resolution of $FM$, the sequence 
$$\dots \to GL_{s}\to GL_{s-1}\to \dots \to GL_{0}\to GM\to 0$$
is a resolution of the $GA$-module $GM$.\\
 \vspace{1em}
\begin{proposition}~\label{homological dimension}
Assume $FA$ is noetherian and complete. 
If $GA$ is of finite (left) global homological dimension, so is $A$. 
\end{proposition}

{\it Proof :} we adjust the proof of [Schn] proposition 10.3.5. to
decreasing filtrations. Let us start by a lemma.

\begin{lemma}\label{Hausdorff and complete}
 If $FN$ is a Hausdorff finite type $FA$-module, then it is complete. 
 \end{lemma}

Indeed, let $FN$ be a finite type $FA$-module. We have a strict exact
sequence 
$$FL={\displaystyle \oplus_{i=1}^{n}FA(-d_{i})}
 \buildrel {p}\over \longrightarrow FN \to 0.$$
The filtration on $FN$ is given by $p(F_{t}L)$. Let us endow the
kernel $K$ of $p$ with the induced topology. We have a strict exact
sequence
$$0 \to FK \to FL \to FN \to 0 .$$
As $N$ is Hausdorff, $K=p^{-1}(\{0\})$ is closed in $FL$. 
The filtered $FA$-module $FN$ is isomorphic to 
$FL/K$, endowed with the quotient topology. Hence, $FN$ is complete
(see lemma 3.0.1).   \\

\begin{lemma}
Assume that $FA$ is noetherian and complete. Then, for any $FA$-module
of finite type $FM$ and any complete $FA$-module $FN$, 
$$\underline{\rm Ext}^{j}_{GA}(GM, GN)=0 \Longrightarrow \Ext^{j}_{A}(M,N)=0.$$ 
\end{lemma}

Let 
$$\cdots \to FL_{n}\to FL_{n-1}\to \cdots \to FL_{0}\to FM \to 0$$
be a filtered resolution of $FN$ by finite free $FA$-modules. Applying
the graduation functor, we get a resolution 
$$\cdots \to GL_{n}\to GL_{n-1}\to \cdots \to GL_{0}\to GM \to 0.$$
Assuming $\underline{\rm Ext}^{j}_{GA}(GM ,GN )=\{0\}$ means that the sequence 
$$\underline{Hom}_{GA}(GL_{j-1}, GN) \to \underline{Hom} _{GA}(GL_{j}, GN) \to 
\underline{Hom} _{GA}(GL_{j+1},GN)$$
is an exact sequence of $GA$-modules. When 
$FL={\displaystyle \mathop \oplus _{i=1}^{n}FA(-d_{i})}$ is finite free, the
$FA$-module $FHom (FL, FN)={\displaystyle \oplus_{i=1}^{n}}FN(d_{i})$
is complete and  
the natural map 
$${\rm G FHom} _{FA}(FL,FN) \to \underline{\Hom}_{GA}(GL ,GN) $$
is an isomorphism. Hence the sequence 
$$\FHom_{FA}(FL_{j-1},FN) \to \FHom_{FA}(FL_{j},FN)\to 
\FHom_{FA}(FL_{j+1},FN)$$
is a strict exact sequence of $FA$-modules (proposition ~\ref{strict exact
  sequences}). When $FL$ is finite free, the underlying module of 
${\rm FHom} _{FA}(FL,FN)$ is $\Hom_{A}(L,N)$. This finishes the proof of the
lemma.  \\

Denote by $d_{GA}$ the (left) global homological dimension of $GA$. 
Let $M$ be a finite type $A$-module. One has an epimorphism 
$$A^{n}\buildrel {p}\over \longrightarrow M \to 0.$$ 
We set $K=\Ker p$. We endow $M$ with the filtration $p(FA^{n})$ and
$K$ with the filtration induced from that of $FA^{n}$. 
We thus define a finite type $FA$-module $FM$ and an $FA$-module $FK$
such that the  exact sequence 
$$0\to FK \to FA^{n }\buildrel {p}\over \longrightarrow FM \to 0$$ 
is strict exact. 
Similarly, we endow $N$ with a filtration $FN$ such that
$FN$ is a finite $FA$-module and we construct a  strict exact sequence 
$$0\to FQ \to FA^{n }\buildrel {q}\over \longrightarrow FN \to 0.$$
The finite $FA$ -module $FN$ is not necessarily Hausdorff but 
$FQ$ is a Hausdorff finitely generated  $FA$-module. 
Hence $FQ$ is a complete $FA$-module
and applying the lemma we get $\Ext_{A}^{j}(M,Q)=0$ when 
$j \geq d_{GA}+1$. 
From the exact sequence 
$$0\to Q \to A^{n }\buildrel {q}\over \longrightarrow N \to 0$$
we get a long exact sequence that shows that 
$\Ext_{A}^{j}(M,N)=0$  when 
$j \geq d_{GA}+1$.  Thus we have showed that : for any finite type $A$
modules $M$ and $N$,  
$$\Ext_{A}^{j}(M,N)=0\;\; {\rm if }\;\; j \geq d_{GA}+1.$$

Let now $N$ be any $A$-module. We have 
$N=\lim\limits_{\rightarrow} N'$ where $N'$ runs over all finitely
generated submodules of $N$. Let $L^{\bullet}$ be a resolution of $M$ by
finitely generated free $A$-modules. We have for all $j \geq d_{GA}+1$
$$\begin{array}{rcl}
\Ext^{j}(M,N)&=& \Ext^{j}(M,\lim\limits_{\rightarrow}N')\\
&=& H^{j}\left ( \Hom_{A}(L^{\bullet},\lim\limits_{\rightarrow}N'
  )\right )\\
&=& H^{j}\left ( \lim\limits_{\rightarrow}\Hom_{A}(L^{\bullet}, N')
\right )\\
&=&  \lim\limits_{\rightarrow}H^{j}\left ( \Hom_{A}(L^{\bullet}, N')
\right )\\
&=& \lim\limits_{\rightarrow} \Ext^{j}_{A}(M, N')
\end{array}$$
where, in the equality before the last equality, we used the fact that
the functor $\lim\limits_{\rightarrow}$ is exact because the set of
finitely generated submodules of $M$ is a directed set 
([Ro] proposition 5.33). Thus we have proved : if $M$ is a finitely generated
$A$-module and $N$ is any $A$-module, then   
$$\Ext^{j}(M,N)=\{0\}\;\;{\rm if }\;\; j\geq d_{GA}+1.$$ 
From this, we deduce ([Ro] theorem 8.16), that the global (left)
dimension of $A$ is finite and inferior or equal to $d_{GA}$. $\Box$

\section{Deformation algebras}

\subsection{Definition and properties}

In this section $k$ will be a field of characteristic $0$ and we will
set $K=k[[h]]$.

\begin{definition}
A  topologically free $K$-algebra $A_{h}$ is a
topologically free $K$-module together with a $K$-bilinear
(multiplication) map $A_{h} \times A_{h} \to A_{h}$ making 
$A_{h}$ into an associative algebra. 

Let $A_{0}$ be an associative $k$-algebra. A deformation of $A_{0}$ is
topologically free $K$-algebra $A_{h}$ such that $A_{0}\simeq A_{h}/hA_{h}$
as algebras. 
\end{definition} 

{\bf Remark :}

If $A_{h}$ is a deformation algebra of $A_{0}$, we may endow it with the
$h$-adic filtration. We then have 
$GA_{h}={\displaystyle \mathop \oplus_{i \in \N}
\frac{h^{i}A_{h}}{h^{i+1}A_{h}}}\simeq A_{0}[h]$ as
$k[h]$-algebra. \\

From proposition 3.0.6, we deduce that 
a deformation algebra of a noetherian algebra is noetherian. \\

{\bf Examples  ([C-P]):}\\

Before giving a list of examples, let us recall the following
definition : 

\begin{definition}
A deformation of a Hopf algebra $(A, \iota , \mu , \epsilon , \Delta
,S)$ over a field $k$ is a topological Hopf algebra 
$(A_{h}, \iota_{h}, \mu_{h}, \epsilon_{h}, \Delta_{h}, S_{h})$ over
the ring $k[[h]]$ such that 

i) $A_{h}$ is isomorphic to $A_{0}[[h]]$ as a $k[[h]]$-module

ii) $A_{h}/hA_{h}$ is isomorphic to $A_{0}$ as Hopf algebra.
\end{definition}

{\it Example 1 : Quantized universal enveloping algebras (QUEA)} 

\begin{definition}

Let $\mathfrak g$ be a Lie bialgebra. A Hopf algebra deformation of 
$U({\mathfrak g})$ , $U_{h}({\goth g})$, such that 
${\displaystyle \frac{U_{h}({\goth g})}{hU_{h}({\goth g})}}$ is
    isomorphic to $U({\goth g})$ as a coPoisson Hopf algebra
 is called a quantization of $U({\goth g})$.

\end{definition}

Quantizations of Lie bialgebras have been constructed
in $[E-K 1]$. \\

{\it Example 2 : Quantization of affine algebraic Poisson  groups}

\begin{definition}
A quantization of an affine algebraic Poisson  group $(G, \{,\})$ is a
Hopf algebra deformation  ${\mathcal F}_{h}(G)$ of the Hopf algebra 
 ${\mathcal F}(G)$ of 
regular functions on $G$, such that 
${\displaystyle \frac {{\mathcal F}_{h}(G)}{h{\mathcal F}_{h}(G)}}$ is isomorphic
to 
$({\mathcal F}(G),\{,\})$ as Poisson Hopf algebra.
\end{definition}

Quantization of affine algebraic Poisson groups have been constructed by
Etingof and Kazhdan ([E-S], 
see also [C-P] for the case where $G$ is simple).\\  

{\it Examples 3: Quantum formal series Hopf algebras (QFSHA)}

The vector space dual $U({\mathfrak g})^{*}$ of the universal
enveloping algebra $U({\goth g})$ of a Lie algebra can be identified
with an algebra of formal power series and it has a natural Hopf
algebra structure, provided we interpret the tensor product 
$U({\mathfrak g})^{*} \otimes U({\mathfrak g})^{*}$ in a suitable
completed sense. If $\goth g$ is a Lie bialgebra, $U({\goth g})^{*}$ is a Hopf Poisson algebra. 

\begin{definition}
A quantum formal series Hopf algebra is a topological Hopf algebra
$B_{h}$ over $k[[h]]$ such that 
${\displaystyle \frac {B_{h}}{hB_{h}}}$ is isomorphic to 
$U({\goth g})^{*}$ as a topological Poisson Hopf algebra for some
finite dimensional Lie bialgebra. 
 
\end{definition}

The following proposition is proved in [K-S] (theorem 2.6)\\
 
\begin{proposition}~\label{flatness}
Let $A_{h}$ be a deformation algebra of $A_{0}$ and let $M$ be an
$A_{h}$-module. Assume that

(i) $M$ has no $h$-torsion

(ii)$M/hM$ is a flat $A_{0}$-module

(iii) $M = {\displaystyle \varprojlim_{n} M/h^{n}M}$

then $M$ is a flat $A_{h}$-module. 
\end{proposition}

%\subsection{Duality functor on deformation algebras}

%We will denote by $A_{h}$ any deformation algebra of a noetherian
%algebra with finite global homological dimension $A_{0}$. 
%Let 
%$Mod \left (A_{h}\right )$ be the abelian category of
%$A_{h}$-modules. 
%Let $D^{b}\left ( Mod (A_{h}) \right )$ be the bounded derived
%category of the abelian category $Mod(A_{h})$. Denote by 
%$D^{b}_{f}\left ( Mod (A_{h}) \right )$ the full subcategory
%of $D^{b}\left ( Mod (A_{h}) \right )$ consisting of
%objects with finitely generated cohomology. 

%Left (respectively right) multiplication endow $A_{h}$ with
%a left (respectively right) $A_{h}$ -module structure. 
%Thus, $A_{h}$ is an 
%$A_{h} \otimes A_{h}^{op}$-module. 
%Introduce the functor $D_{A_{h}}$ from 
%$D^{b}_{f}\left (  Mod \left (A_{h} \right ) \right )$ 
%to $D^{b}_{f}\left (  Mod \left (A_{h}^{op} \right ) \right )$
%$$\forall M^{\bullet} \in D^{b}_{f}\left (  A_{h} \right ),\;\; 
%D_{A_{h}}(M^{\bullet})=
%RHom_{A_{h}}\left ( M^{\bullet} , A_{h} \right ). $$
%The canonical arrow 
%$M^{\bullet} \to D_{A^{op}_{h}} \circ D_{A_{h}}(M^{\bullet})$ 
%is an isomorphism. To show this, prove it for $M=A_{h}$ and
%use standard homological algebra arguments.\\

\section{A quantization of the character $trad$ }

\begin{theorem} \label{qtrad}
Let $A_{0}$ be a noetherian $k$-algebra 
%with finite global homological dimension 
and let $A_{h}$ be a
deformation of $A_{0}$. 
Assume that $k$ has a left 
$A_{0}$-module structure such that   there exists an integer $d$ such
that 
$$\left \{\begin{array}{l}
Ext^{i}_{A_{0}}\left (k, A_{0} \right )=\{0\}\;{\rm if}\; i\neq d \\
Ext^{d}_{A_{0}}\left (k, A_{0} \right )\simeq k
\end{array}\right .$$
Assume that $K$ is endowed with a $A_{h}$-module structure which
reduces modulo $h$ to the $A_{0}$-module structure on $k$ we started with. 
Then 

a) $Ext^{i}_{A_{h}}\left (K, A_{h} \right )$ is
zero if $i \neq d$.

b) $Ext^{d}_{A_{h}}\left ( K, A_{h} \right )$ is
a free $K$-module of dimension 1. 
By right multiplication, it is a right $A_{h}$-module. 
It is a lift of the right $A_{0}$-module structure (given by right
multiplication ) on $Ext^{d}_{A_{0}}\left (k, A_{0} \right )$.
\end{theorem}

{\bf Notation :} The right $A_{h}$-module  
$Ext^{d}_{A_{h}}\left (k, A_{h} \right )$ will be denoted
$\Omega_{A_{h}}$ and the character defined by this action $\theta_{A_{h}}$. \\

{\bf Remark :} In [K-S] (paragraph 6), Kashiwara and Schapira make a similar construction in the set up of $DQ$-algebroids. In [C2], it is shown that a result similar to theorem ~\ref{qtrad} holds for 
$U_{q}({\goth g})$ (${\goth g}$ semi-simple).\\

{\it Example 1 : Quantized universal enveloping algebras}

Poincar{\'e} duality gives us the following result for any finite
dimensional Lie algebra. 

$$\left \{ \begin{array}{l}
Ext^{i}_{U({\mathfrak g})}\left (k, U(\mathfrak g) \right ) =\{0\} 
\;\;{\rm  if} i\neq 0\\
Ext^{dim \goth g}_{U({\mathfrak g})}\left (k, U(\mathfrak g) \right ) \simeq
\Lambda^{dim {\mathfrak g}}({\mathfrak g} ^{*})\;\;. 
\end{array}\right .$$
The character defined by the right action of $U({\goth g})$ on 
$Ext^{dim \goth g}_{U({\mathfrak g})}\left (k, U(\mathfrak g) \right
)$ is $trad_{\mathfrak g}$ ([C1]). Thus, 
the character defined by the theorem ~\ref{qtrad} is a quantization of 
the character $trad_{\mathfrak g}$. 

$\bullet$ If ${\mathfrak g}$ is a complex semi-simple algebra, as 
$H^{1}({\mathfrak g}, k)=\{0\}$ ([H-S] p 247), there exists a unique
lift of the trivial representation of $U_{h}({\mathfrak g})$, hence 
the representation $\Omega_{U_{h}({\mathfrak g})}$ is the trivial representation.  

$\bullet$ Let $\mathfrak a$ be a $k$-Lie algebra. Denote
by ${\mathfrak a}_{h}$ the Lie algebra obtained from $\mathfrak a$ by
multiplying the bracket of ${\mathfrak  a}$ by $h$. 
Thus, for any elements $X$ and $Y$ of 
${\mathfrak a}_{h}\simeq {\mathfrak a}$, 
$$[X,Y]_{{\mathfrak a}_{h}}=h[X,Y]_{\mathfrak a}.$$
Denote by 
$\widehat {U({\mathfrak a}_{h})}$ the $h$-adic completion of 
$U({\mathfrak a}_{h})$. Then $\widehat {U({\mathfrak a}_{h})}$ is a
Hopf deformation of $( {\mathfrak a}^{ab}, \delta =0 )$. The character
$\theta_{\widehat {U({\mathfrak a}_{h})}}$ 
defined by the theorem in this case is given by 
$$\forall X \in {\mathfrak a}, \;\; 
\theta_{\widehat {U({\mathfrak a}_{h})}}(X)=htrad_{\mathfrak a}(X).$$  
Thus, even if ${\mathfrak g}$ is
unimodular, the character defined by the right action of
$U_{h}({\mathfrak g})$ on 
$\Omega_{U_{h}({\mathfrak g})}
\simeq \wedge ^{dim{\mathfrak g}}({\mathfrak g}^{*})[[h]]$
might not be trivial. 

$\bullet$ We consider the following Lie algebra : 
 ${\goth a}={\displaystyle \mathop
  \oplus_{i=1}^{5}ke_{i}}$ with non zero bracket 
$\left [e_{2},e_{4} \right ]=e_{1}$.  
Consider  $k[[h]]$-Lie algebra structure on ${\goth a }[[h]]$ 
defined by the following non zero brackets 
$$\begin{array}{l}
\left [e_{3},e_{5}\right ]=he_{3}\\
\left [ e_{2},e_{4} \right ]=2e_{1}
\end{array}$$
$\widehat{U({\goth a}[[h]])}$ is a quantization of $U({\goth a})$. It
is easy to see that  
$$\begin{array}{l}
\theta_{\widehat{U({\goth a}[[h]])}}(e_{i})=0\;\;{\rm if}\;\;i\neq 5\\
\theta_{\widehat{U({\goth a}[[h]])}}(e_{5})=-h.
\end{array}$$

{\it Example 2} 

The theorem  ~\ref{qtrad} also applies to quantization of affine algebraic
Poisson  groups. If $G$
is an affine algebraic Poisson group with neutral element $e$,
 we take $k$ to be given by the counit of the Hopf
algebra ${\mathcal F}(G)$. One has [A-K] 
$$\begin{array}{l}
Ext^{i}_{{\mathcal F}(G)}\left ( k, {\mathcal F}(G) \right )=\{0\}\;{\rm
  if}\; i\neq dim G \\
Ext^{dim G}_{{\mathcal F}(G)}\left ( k, {\mathcal F}(G) \right )\simeq 
 \wedge ^{dim G}\left (
\left ({\mathcal M}_{e}/{\mathcal M}_{e}^{2} \right )^{*} \right )
\end{array}$$
where
$${\mathcal M}_{e}=\{f \in {\mathcal F}(G) \mid f(e)=0\}.$$
Let $\mathfrak g$ be a real Lie algebra.  
The algebra of regular functions on ${\mathfrak g}^{*}$, 
${\mathcal F}({\mathfrak g}^{*})$,  is isomorphic to
$S({\mathfrak g})$ and is naturally equipped with a Poisson structure
given by :
$$\forall X,Y \in {\mathfrak g}, \{X,Y\}=[X,Y].$$  
In the example above, $\widehat {U({\mathfrak g}_{h})}$ is a
quantization of the Poisson algebra ${\mathcal F}({\mathfrak g}^{*})$.  
${\mathcal F}({\mathfrak g}^{*})$ acts trivially on 
$Ext^{dim {\mathfrak g}}_{{\mathcal F}({\mathfrak g}^{*})}
\left ( k, {\mathcal F}({\mathfrak g}^{*}) \right )$
whereas the action of 
${\mathcal F}_{h}({\mathfrak g}^{*})\simeq\widehat {U({\mathfrak
    g}_{h})}$ on 
$Ext^{dim {\mathfrak g}}_{{\mathcal F}_{h}({\mathfrak g}^{*})}
\left ( k, {\mathcal F}_{h}({\mathfrak g}^{*}) \right )$ is not trivial. \\

{\it Example 3 :} The theorem  ~\ref{qtrad} also applies to quantum formal series Hopf algebras. \\

{\it Proof of the theorem ~\ref{qtrad}:}

Let us consider a resolution of the $A_{h}$-module $K$ by
filtered finite free $A_{h}$-modules 
$$\begin{array}{l}
\dots \to FL^{i+1}\buildrel {\partial_{i+1}} \over \longrightarrow FL^{i}
\buildrel{\partial_{i}}\over \longrightarrow  \dots \buildrel{\partial_{2}} \over
\longrightarrow   FL^{1} \buildrel{\partial_{1}} \over 
\longrightarrow FL^{0}\to K \to \{0\}\\
FL^{i}={\displaystyle \oplus _{k=1}^{d_{i}}FA_{h}(-m_{j,i})}
\end{array}$$
so  that the graded complex 
$$\dots GL^{i+1} \buildrel {G\partial_{i+1}} \over \longrightarrow 
GL^{i}\buildrel {G\partial_{i}} \over \longrightarrow 
 \dots \to GL^{1}\buildrel {G\partial_{1}} \over \longrightarrow 
  GL^{0}\to k[h] \to \{0\}$$
is a resolution of the $A_{0}[h]$-module $k[h]$.
Consider the complex 
$M^{\bullet}=\left ( Hom_{A_{h}}(L^{\bullet}, A_{h}), ^{t}\partial^{\bullet} \right )$.
Recall that  there is a natural filtration on 
 $Hom_{A_{h}}(L^{i}, A_{h})$ defined by 
$$F_{t}Hom_{A_{h}}(L^{i}, A_{h})=
\{\lambda \in Hom_{A_{h}}(L^{i}, A_{h})
\mid \lambda \left ( F_{p}L^{i} \right )\subset
F_{t+p}A_{h}\}.$$
 One has 
an isomorphism of right $FA$-modules 
$${\rm FHom}_{A_{h}}(L^{i}, A_{h})=
{\displaystyle \mathop \oplus_{j=1}^{r_{i}}FA(m_{j,i})} $$
Hence 
$${\rm GFHom}_{A_{h}}(L^{i}, A_{h})
\simeq \underline{  Hom}_{GA_{h}}(GL^{i}, GA_{h})$$ and the
complex $\underline{\rm Hom}_{GA_{h}}(GL^{i}, GA_{h})$
computes $\underline{Ext}^{i}_{GA_{h}}\left (k[h], GA_{h} \right )$. We have the
following isomorphisms of right $A_{0}[h]$-modules. 
$$\underline{\rm Ext}^{i}_{GA_{h}}\left (k[h], GA_{h} \right )\simeq 
\underline{\rm Ext}^{i}_{A_{0}[h]}\left (k[h], A_{0}[h] \right )\simeq 
\Ext^{i}_{A_{0}}(k,A_{0})[h].$$

If $i \neq d$, then 
$\underline{\Ext}^{i}_{GA_{h}}\left (k[h], GA_{h} \right )
=\{0\}$. This means that the sequence 
$$\underline{Hom}_{GA}(GL_{i-1}, GA_{h}) 
\buildrel{^{t}G\partial_{i}} \over \longrightarrow
 \underline{Hom} _{GA}(GL_{i}, GA_{h}) 
  \buildrel{^{t}G\partial_{i+1}} \over \longrightarrow
\underline{Hom} _{GA}(GL_{j+1},GA_{h})$$
is an exact sequence of $GA_{h}$-modules.  
Hence, applying \ref{strict exact sequences} the sequence 
$$\FHom_{FA}(FL_{i-1},FN) 
\buildrel{^{t}\partial_{i}} \over \longrightarrow
 \FHom_{FA}(FL_{i},FN)
\buildrel{^{t}\partial_{i+1}}  \over \longrightarrow
\FHom_{FA}(FL_{i+1},FN)$$ 
is strict exact. 
As  $FL_{i}$ is finite free, the underlying module of 
${\rm FHom} _{FA}(FL_{i},FN)$ is $\Hom_{A}(L_{i},N)$.  Hence we have proved that 
$Ext^{i}_{A_{h}}(K,A_{h})=\{0\}$ if $i \neq d$. 

We have also proved that all the maps 
$^{t}\partial_{i}$ are strict. 
Hence, by proposition ~\ref{spectral sequence}, we have for all integer $i$ 
$$GExt^{i}_{A_{h}}(k[[h]],A_{h}) \simeq \underline{Ext}^{i}_{GA_{h}}(k[h], A_{0}[h])
\simeq Ext ^{i}_{A_{0}}(k,A_{0})[h]$$

As $Im ^{t}\partial_{i-1}$, endowed with the induced topology, 
 is a finite type Hausdorff $FA_{h}$-module, it is a complete $FA_{h}$-module
 (see lemma  \ref{Hausdorff and complete}). Hence it is closed 
 in $Ker ^{t}\partial_{i}$ and the $\Ext^{i}_{A_{h}}\left (K, A_{h} \right )$'s are Hausdorff. 

As $\Ext^{d}_{A_{h}}\left (K, A_{h} \right )$ is  Hausdorff and 
$GExt^{d}_{A_{h}}(k[[h]],A_{h}) \simeq Ext ^{d}_{A_{0}}(k,A_{0})[h]$, the $k[[h]]$-module 
$\Ext^{d}_{A_{h}}\left (K, A_{h} \right )$ is a one dimensional. 

This finishes the proof of the theorem ~\ref{qtrad}. $\Box$\\

%Let us now compute the 
%$Ext^{i}\left ( {\displaystyle \frac{K}{h^{p}K}}, U_{h}({\goth
%    g})\right )$.
%One has the exact sequence 
%$$0\to K \buildrel{\mu_{h^{p}}}\over \longrightarrow K \to 
%{\displaystyle \frac{K}{h^{p}K}} \to 0$$
%from which we deduce the following long exact sequence 
%$$\begin{array}{l}
%0 \to Ext^{0}_{U({\goth g})}\left ({\displaystyle \frac{K}{h^{p}K}},
% U_{h}({\goth g}) \right ) \to  
%Ext^{0}_{U({\goth g})}\left ( K, U_{h}({\goth g}) \right ) \to \dots\\ 
%\dots \to Ext^{i}_{U({\goth g})}\left ({\displaystyle \frac{K}{h^{p}K}},
%  U_{h}({\goth g}) \right ) \to  
%Ext^{i}_{U({\goth g})}\left (K, U_{h}({\goth g}) \right ) \to 
%Ext^{i}_{U({\goth g})}\left (K, U_{h}({\goth g}) \right ) \\
% \to Ext^{i+1}_{U({\goth g})}\left ({\displaystyle \frac{K}{h^{p}K}},
%  U_{h}({\goth g}) \right ) \dots \\
%\end{array}$$
%From this long exact sequence, we deduce that 
%\begin{itemize}
%\item
%$Ext^{i}_{U({\goth g})}
%\left (\frac{K}{h^{p}K}, U_{h}({\goth g}) \right )\neq 0
%\;\;{\rm if} \;\; i = dim{\goth g}, dim{\goth g}+1.$ 
%\item
%$Ext^{dim \goth g}_{U({\goth g})}
%\left (\frac{K}{h^{p}K}, U_{h}({\goth g}) \right )$
%and 
%$Ext^{dim {\goth g}+1}_{U({\goth g})}\left (\frac{K}{h^{p}K}, U_{h}({\goth g})
%\right )$ are torsion modules. 
%\end{itemize}
%This is in contradiction  with the fact that 
%$D^{2}_{U_{h}({\goth  g})}(K)\simeq K$. 
%Thus $Ext^{dim \goth g}_{U({\goth g})}\left (K, U_{h}({\goth g})
%\right )$ is isomorphic to $K$ and the proposition is proved. 

{\it From now on, we assume that $A_{h}$ is a topological Hopf algebra
  and that its action on $K$ is given by the counit. The antipode of $A_{h}$ will be denoted 
  $S_{h}$.}\\

If $V$ is a left $A_{h}$-module, we set $V^{r}$ (respectively $V^{\rho}$) the right $A_{h}$-module defined by 
$$\forall a \in A_{h}, \; \forall v \in V,\;
v \cdot_{S_{h}} a = S_{h}(a) \cdot v \;\; 
({\rm respectively } \;\; v \cdot_{S_{h}^{-1}} a = S_{h}^{-1}(a) \cdot v) .$$ 
Similarly, if $W$ is a right $A_{h}$-module, we set $W^{l}$ 
(respectively $W^{\lambda}$) the left  $A_{h}$-module defined by 
$$\forall a \in A_{h}, \; \forall w \in W,\;
a \cdot_{S_{h}} w = w \cdot S_{h}(a) \;\; 
({\rm respectively } \;\; a \cdot_{S^{-1}_{h}} w = w \cdot S_{h}^{-1}(a) ) .$$ 
One has $\left (V^{r}\right )^{\lambda}=V$,
$\left (V^{\rho}\right )^{l}=V$, 
  $\left ( W^{l}\right )^{\rho}=W$ and  $\left ( W^{\lambda}\right )^{r}=W$. 
Thus, we have   
defined two (in the case where $S_{h}^{2}\neq id $)  equivalences of categories between the category of left  
$A_{h}$-modules and the category of right $A_{h}$-modules, that is to
say left  $A_{h}^{op}$-modules. \\

Let $Mod \left (A_{h}\right )$ be the abelian category of left 
$A_{h}$-modules and  $D\left ( Mod (A_{h}) \right )$ be the  derived
category of the abelian category $Mod(A_{h})$.  We may consider 
 $A_{h}$ as an 
$A_{h} \otimes A_{h}^{op}$-module. 
Introduce the functor $D_{A_{h}}$ from 
$D\left (  Mod \left (A_{h} \right ) \right )$ 
to $D\left (  Mod \left (A_{h}^{op} \right ) \right )$
$$\forall M^{\bullet} \in D\left (  A_{h} \right ),\;\; 
D_{A_{h}}(M^{\bullet})=
RHom_{A_{h}}\left ( M^{\bullet} , A_{h} \right ). $$
If $M$ is a finitely generated module,  the canonical arrow 
$M \to D_{A^{op}_{h}} \circ D_{A_{h}}(M)$ 
is an isomorphism. \\

Let $V$ be a left $A_{h}$-module, then, by transposition, 
 $V^{*}=Hom_{K}(V,K)$ is naturally endowed with a  right 
$A_{h}$-module structure. Using the antipode, we can also see it as a
left module structure. Thus, one has : 
$$\forall u \in A_{h}\, \forall f \in V^{*}, \;\; u \cdot f =f \cdot S_{h}(u).$$
We endow $ \Omega_{A_{h}}\otimes V^{*}$ with the following
right  $A_{h}$-module structure :
$$\begin{array}{l}
\forall u \in A_{h}\, \forall f \in V^{*}, \;
\forall \omega \in \Omega_{A_{h}},\\
(\omega \otimes w)\cdot u =
 \lim\limits_{n\to +\infty}
{\displaystyle \mathop \sum_{j}}
\theta_{A_{h}} (u'_{j,n})\omega  \otimes 
f \cdot S_{h}^{2}(u''_{j,n}) 
\end{array}$$
where $\Delta (u)=\lim\limits_{n\to +\infty}
{\displaystyle \mathop \sum_{j}u'_{j,n}\otimes u''_{j,n}}$.\\

\begin{theorem} \label{two dualities}
Let $V$ be an $A_{h}$-module free of finite type as a 
$k[[h]]$-module. Then  
$D_{A_{h}}(V)$ and $V^{*}\otimes \Omega_{A_{h}}$ are
isomorphic in $D\left ( A_{h}^{op} \right )$. 
\end{theorem}

{\it Proof of the theorem :}

In the proof of this theorem, we will make use of the following lemma 
(see [Du1], [C1]).

\begin{lemma}\label{bimodule isomorphism}
Let $W$ be a left $A_{h}$-module. 
$ A_{h}\widehat{\otimes} W$ is endowed with two different
structures of 
$A_{h}\otimes A_{h}^{op}$-modules. 
The first one denoted 
$\left ( A_{h} \widehat{\otimes} W\right )_{1}$ is described as
follows : Let $w$ be an element of $W$ and let $u,a$ be two
elements of $A_{h}$. We set 
$\Delta (a)=
\lim\limits_{n \to +\infty}{\displaystyle \sum_{i}}a'_{i,n}\otimes a''_{i,n}$.
Then 
$$\begin{array}{l}
(u \otimes w)\cdot a= ua \otimes w\\
a \cdot (u \otimes w)=\lim\limits_{n \to + \infty}
{\displaystyle \sum_{i}}a'_{i,n}u \otimes a''_{i,n}\cdot w
\end{array}$$
The second one   denoted 
$\left ( A_{h}\widehat{\otimes} W\right )_{2}$ is described as
follows : 
Then 
$$\begin{array}{l}
a \cdot (u \otimes w)=au\otimes w\\
(u \otimes w)\cdot a=\lim\limits_{n \to + \infty}
{\displaystyle \sum_{i}}ua'_{i,n}\otimes  S_{h}(v''_{i,n})\cdot w
\end{array}$$ 
The $A_{h}\otimes A_{h}^{op}$-modules
$\left ( A_{h}\widehat{\otimes} W \right )_{1}$ and 
$\left ( A_{h}\widehat{\otimes} W \right )_{2}$ are isomorphic.
\end{lemma}
{\it Proof of the lemma :}

The map 
$$\begin{array}{rcl}
\Psi : \left ( A_{h}\widehat{\otimes} W \right )_{2} & \to & 
\left ( A_{h}\widehat{\otimes} W \right )_{1}\\
 u \otimes w& \mapsto &
\lim\limits_{n\to +\infty}{\displaystyle \sum_{i}}u'_{i,n} \otimes u''_{i,n}\cdot w
\end{array}$$
where 
$\Delta (u)=
\lim\limits_{n \to +\infty}{\displaystyle \sum_{i}}u'_{i,n}\otimes
u''_{i,n}$ is an isomorphism 
of $A_{h}\otimes A_{h}^{op}$-modules from 
$\left ( A_{h}\widehat{\otimes}W \right )_{2}$ to 
$\left (  A_{h}\widehat{\otimes}W \right )_{1}$.  Moreover 
$$\Psi ^{-1}(u \otimes w)= \sum u'_{i,n} \otimes S_{h}(u''_{i,n})\cdot w.$$
This finishes the proof of the lemma. \\

Let $L^{\bullet}$ be a resolution of $K$ by free 
$A_{h}$-modules. We endow $L^{i}\otimes V$ with the following left $A_{h}$-module structure :
$$a \cdot (l \otimes v)=
\lim\limits_{n\to +\infty}{\displaystyle \sum_{i}}a'_{i,n}\cdot l  \otimes a''_{i,n}\cdot v. $$
Then $L^{\bullet}\otimes V$ is a
resolution of $V$ by free $A_{h}$-modules. Using the relation 
$$a \cdot l \otimes v=\lim\limits_{n\to +\infty}{\displaystyle \sum_{i}}a'_{i,n}\left [
l \otimes S_{h}(a''_{i,n})\cdot v
\right ] $$
 one shows the
following sequence of $A_{h}$-isomorphisms 
$$\begin{array}{rcl}  
$$D_{A_{h}}(V) & \simeq &
Hom_{A_{h}}\left (L \otimes V, A_{h} \right )\\
& \simeq & Hom _{A_{h}} \left ( L, 
\left ( A_{h}\otimes V^{*}\right )_{1} \right )\\
& \simeq &
Hom _{A_{h}} \left ( L, 
\left ( A_{h}\otimes V^{*}\right )_{2} \right )\\
&\simeq &
 RHom_{A_{h}}(K, A_{h}) \otimes V^{*}. \Box
\end{array}$$

\section{Link with quantum duality}

\subsection{Recollection on the quantum dual principle.}

The quantum dual principle ([Dr], see [G] for a detailed treatment)
states that there exist two functors, namely 
$\left ( \right )': QUEA \to QFSA$ and 
$\left ( \right )^{\vee}: QFSA \to  QUEA$ which are
inverse of each other. If $U_{h}({\goth g})$ is a quantization of
$U({\goth g})$ and $F_{h}[[{\goth g}]]$ is a quantization of 
$F[[{\goth g}]]=U({\goth g})^{*}$, then $U_{h}({\goth g})'$ is a
quantization of $F[[{\goth g}^{*}]]$ and 
$F_{h}[[{\goth g}]]^{\vee}$ is a quantization of $U({\goth g}^{*})$. 

Let's recall the construction of the functor 
$\left ( \right )^{\vee}: QFSA \to QUEA$ 
which is the one we will need.  Let ${\goth g}$ be a Lie bialgebra and 
$F_{h}[[{\goth g}]]$ a quantization of $F[[{\goth g}]]=U({\goth
  g})^{*}$. 
For simplicity we will write $F_{h}$ instead of $F_{h}[[{\goth g}]]$.
If  $\epsilon _{h}$ denotes 
the counit of $F_{h}$, set  
$I:=\epsilon_{h}^{-1}(hk[[h]])$ and $J=Ker \epsilon_{h}$. Let 
$$F_{h}^{\times}:=
{\displaystyle \sum_{n \geq 0}}h^{-n}I^{n}=
 {\displaystyle \sum_{n \geq 0}}(h^{-1}I)^{n}
={\bigcup _{n \geq 0}}(h^{-1}I)^{n}$$
be the  $k[[h]]$-subalgebra of 
$k((h)){\displaystyle \mathop \otimes_{k[[h]]}}F_{h}$ generated by $h^{-1}I$. 
As $I=J+hF_{h}$, one has 
$F_{h}^{\times}={\displaystyle \sum_{n \geq 0}}h^{-n}J^{n}$. 
Define $F_{h}^{\vee}$ to be the $h$-adic completion of
the $k[[h]]$-module $F_{h}^{\times}$. The coproduct (respectively
counit,  antipode) on $F_{h}$
provides a coproduct (respectively counit, antipode) 
on $F_{h}^{\vee}$ and $F_{h}^{\vee}$ is endowed with  a Hopf
algebra structure. A precise
description of  $F_{h}^{\vee}$ is given in [G]. Let us
recall it as we will need it for our computations.  The 
algebras $F_{h}/hF_{h}$ and $k[[\bar{x}_{1},\dots,  \bar{x}_{n}]]$ 
are isomorphic. We denote  $\pi : F_{h} \to F_{h}/hF_{h}$ be the natural
projection. We may choose $x_{j} \in \pi^{-1}(\bar{x}_{j})$ for any
$j$ such that $\epsilon_{h}(x_{j})=0$,
then $F_{h}$ and $k[[x_{1},\dots,x_{n},h]]$ are isomorphic as $k[[h]]$-
topological module and $J$ is the set of formal series $f$ whose
degree in the $x_{j}$, $\partial_{X}(f)$ (that is the degree of the
lowest degree monomials occuring in the series with non zero
coefficients) is strictly positive. 
As $F_{h}/hF_{h}$ is commutative, one has 
$$x_{i}x_{j}-x_{j}x_{i}=h\chi_{i,j}$$
with $\chi_{i,j} \in F_{h}$. As $\chi_{i,j}$ is in $J$, it can be
written as follows : 
$$\chi_{i,j}={\displaystyle \sum_{a=1}^{n}}c_{a}(h)x_{a}+
f_{i,j}(x_{1}, \dots ,x_{n},h)$$
with $\partial_{X}(f_{i,j})>1$.  
If $\check{x}_{i}=h^{-1}x_{j}$, then 
$$F_{h}^{\vee}=\{
f={\displaystyle \mathop \sum_{r \in \N}} P_{r}(\check{x}_{1},\dots, 
\check{x}_{n})h^{r}\mid P_{r}(X_{1},\dots ,X_{n}) \in
k[X_{1},\dots,X_{n}]\}. $$
Thus $F_{h}^{\vee}$ and $k[\check{x}_{1},\dots, \check {x}_{n}][[h]]$
are isomorphic as a topological $k[[h]]$-modules. One has   
$$\check{x}_{i}\check{x}_{j}-\check{x}_{j}\check{x}_{i}=
{\displaystyle \mathop \sum_{a=1}^{n}}c_{a}(h)\check{x}_{a}+ h^{-1}
\check{f}_{i,j}(\check{x}_{1}, \dots, \check{x}_{n},h)$$
where $\check{f}_{i,j}(\check{x}_{1}, \dots, \check{x}_{n},h)$ is 
obtained from $f_{i,j}(x_{1}, \dots , x_{n})$ by writing 
$x_{j}=h\check{x}_{j}$. The element 
$h^{-1}\check{f}_{i,j}(\check{x}_{1}, \dots, \check{x}_{n},h)$ is in 
$hk[\check{x}_{1},\dots, \check{x}_{n}][[h]]$. 
The $k$-span of the set of cosets 
$\{e_{i}=\check{x}_{i} \;{\rm  mod}\; hF_{h}^{\vee} \}$ is a Lie algebra
isomorphic to ${\goth g}^{*}$. The map
$\Psi : F_{h}^{\vee} \to  U({\goth g}^{*})[[h]]$ defined by 
$$\Psi \left ( {\displaystyle \mathop \sum_{r \in \N}}
P_{r}(\check{x} _{1}, \dots ,\check{x}_{n})h^{r}  \right )= 
{\displaystyle \mathop \sum_{r \in \N}}
P_{r}(e_{1}, \dots ,e_{n})h^{r}$$
is an isomorphism of topological $k[[h]]$-modules. The algebra 
${\displaystyle \mathop \frac {F_{h}^{\vee}}{hF_{h}^{\vee}}}$ is
  isomorphic to $U({\mathfrak g}^{*})$ and $F_{h}({\mathfrak g})^{\vee}$ is a
  quantization of the coPoisson Hopf algebra $U({\mathfrak g}^{*})$.   
Denote by $\cdot_{h}$ multiplication on $F_{h}$ and its transposition
to $U({\goth g}^{*})[[h]]$ by $\Psi$. To compute 
$e_{1}^{a_{1}}\dots e_{n}^{a_{n}} \cdot_{h} e_{1}^{b_{1}}\dots
e_{n}^{b_{n}}$ we proceed as follows : we compute 
$\check{x}_{1}^{a_{1}}\dots \check{x}_{n}^{a_{n}} \cdot_{h} 
\check{x}_{1}^{b_{1}}\dots \check{x}_{n}^{b_{n}}$ in $F_{h}^{\vee}$ 
and write it under  the
form ${\displaystyle \mathop \sum_{r \in \N}}
P_{r}(\check{x} _{1}, \dots ,\check{x}_{n})h^{r}$. Then 
$$e_{1}^{a_{1}}\dots e_{n}^{a_{n}} \cdot_{h} e_{1}^{b_{1}}\dots
e_{n}^{b_{n}}=
{\displaystyle \mathop \sum_{r \in \N}}
P_{r}(e_{1}, \dots ,e_{n})h^{r}.$$
If $u$ and $v$ are in $U({\goth g}^{*})[[h]]$, one writes 
$u*v={\displaystyle \mathop \sum_{r \in \N}}h^{r}\mu_{r}(u,v)$. One
knows that the first non zero $\mu_{r}$ is a 1-cocycle of the Hochschild
cohomology. 

If $P$ in $k[X_{1}, \dots, X_{n}]$ can be written 
$P={\displaystyle \mathop \sum _{i_{1},\dots ,i_{n}}}a_{i_{1},\dots, i_{n}}
X_{1}^{i_{1}}\dots X_{n}^{i_{n}}$, one sets 
$$P^{\otimes}(e_{1},\dots, e_{n})= 
{\displaystyle \mathop \sum _{i_{1},\dots ,i_{n}}}a_{i_{1},\dots, i_{n}}
e_{1}^{\otimes i_{1}}\dots e_{n}^{\otimes i_{n}}$$
and if $g \in k[X_{1},\dots ,X_{n}][[h]]$ can be written 
$g={\displaystyle \mathop \sum _{i=1}^{r}}P_{r}(X_{1},\dots,
X_{r})h^{r}$,
then one sets :
$$g^{\otimes}(e_{1},\dots ,e_{n})=
{\displaystyle \mathop \sum _{i=1}^{r}}P_{r}^{\otimes}(e_{1},\dots,
e_{r})h^{r} .$$
{\it Fact :}

{\it $(F_{h})^{\vee}$ is isomorphic as an algebra to 
$$U_{h}({\goth g}^{*})\simeq 
\frac{T_{k[[h]]}\left ( 
{\displaystyle \mathop \oplus_{i=1}^{n}} k[[h]]
e_{i} \right )}{I}$$
where $I$ is the closure (in the $h$-adic topology ) of the two sided
ideal generated by the relations 
$$e_{i}\otimes e_{j}-e_{j}\otimes e_{i}= 
{\displaystyle \mathop \sum_{k=1}^{n}}c_{k}(h)e_{k}+
h^{-1}\check{f}^{\otimes}_{i,j}(e_{1}, \dots, e_{n},h).$$}
 Let us prove this fact. 
Let $\Omega : T_{k[[h]]}
\left ( {\displaystyle \mathop \oplus_{i=1}^{n}}k[[h]]e_{i} \right ) 
 \to  F_{h}^{\vee}$ that sends $e_{i}$ to 
$\check{x}_{i}$.
One has $I\subset Ker\Omega $  and we need to prove that 
$Ker \Omega \subset I$. Let ${\mathcal R}$ be in \linebreak  
$T_{k[[h]]}\left ( {\displaystyle \mathop \oplus_{i=1}^{n}}
  k[[h]]e_{i} \right ) $  
be such that 
$\Omega ({\mathcal R})=0$. Then, modulo $h$, we get 
$\overline{\Omega}(\overline{\mathcal R})=0$. Hence there exist 
$(u_{i,j}^{0})$ and $(v_{i,j}^{0})$ in 
$T_{k}({\displaystyle \mathop \oplus_{i=1}^{n}ke_{i}})$ such that 
$$\overline{\mathcal R}=
{\displaystyle \sum_{i,j}}u_{i,j}^{0}\otimes \left ( e_{i}\otimes
  e_{j}-e_{j}\otimes e_{i}-[e_{i},e_{j}] \right ) \otimes v^{0}_{i,j}$$
and  ${\mathcal R}-
{\displaystyle \sum_{i,j}}u_{i,j}^{0}\otimes \left ( 
e_{i}\otimes e_{j}-e_{j}\otimes e_{i}
-{\displaystyle \sum_{a=1}^{n}}c_{a}(h)e_{a}
-h^{-1}\check{f}_{i,j}^{\otimes}(e_{1},\dots , e_{n}, h) 
\right ) \otimes v_{i,j}^{0} \in  hKer\Omega$. 
Hence there exist ${\mathcal R}_{1} \in Ker \Omega$ be such that 
$${\mathcal R}-
{\displaystyle \sum_{i,j}}u_{i,j}^{0}\otimes \left ( 
e_{i}\otimes  e_{j}-e_{j}\otimes e_{i}-
{\displaystyle \sum_{a=1}^{n}}c_{a}(h)e_{a}- 
h^{-1}\check{f}_{i,j}^{\otimes}(e_{1},\dots ,e_{n},h ) \right ) 
\otimes v^{0}_{i,j} =h{\mathcal R}_{1}.$$ 
Reproducing the same reasoning, we find 
$\left ( u_{i,j}^{1} \right )$ and 
$\left ( v_{i,j}^{1} \right ) $ in 
$T_{k}\left 
( {\displaystyle \mathop \oplus_{i=1}^{n} ke_{i}} \right )$  
such that 
$${\mathcal R}_{1}-
{\displaystyle \sum_{i,j}}u_{i,j}^{1}\otimes 
\left ( e_{i}\otimes  e_{j}-e_{j}\otimes e_{i}
-{\displaystyle \sum_{a=1}^{n}}c_{a}(h)e_{a}
-h^{-1}\check{f}_{i,j}(e_{1},\dots, e_{n},h) 
\right ) \otimes v_{i,j}^{1} =h{\mathcal
  R}_{2}$$
and going on like this, we show that ${\mathcal R}$ is in 
$I$. \\

%Let ${\mathcal M}$ be the maximal ideal of the local algebra $F[[{\goth
%  g}]]$. Let $V[[h]]$ be an $F_{h}$-module such that ${\mathcal M}$ acts
%trivially on the $F[[{\goth g}]]$-module $V$. If $\alpha$ is in $J$
%and $v \in V[[h]]$, then $\alpha \cdot v$ is in $hV[[h]]$. Thus
%$V[[h]]$ can be considered as a 
%$F_{h}^{\vee}\simeq U_{h}({\goth g}^{*})$-module by
%putting :
%$$\forall v \in V, \;\;\check{x}_{i}\cdot v=h^{-1}x_{i}\cdot v .$$
%Considered as an $U_{h}({\goth g}^{*})$, $V[[h]]$ will be denoted 
%$V[[h]]^{\vee}$. 

\subsection{Deformation of the Koszul complex}

Let ${\goth a}$ be a $k$-Lie algebra. There is  a well known
resolution of the trivial $U({\goth a})$-module, namely the Koszul
resolution 
$K=\left ( U({\goth a})\otimes \wedge^{\bullet}{\goth a}, \partial \right )$
where 
$$\begin{array}{l}
\partial \left (u \otimes X_{1}\wedge \dots \wedge X_{n}\right )=
{\displaystyle \mathop \sum_{i=1}^{n}}(-1)^{i-1}uX_{i} \otimes 
X_{1}\wedge \dots \wedge \widehat{X_{i}}\wedge \dots \wedge X_{n}\\
\;\;\;\;{\displaystyle \sum_{i<j}(-1)^{i+j}}
u \otimes 
[X_{i},X_{j}]\wedge X_{1}\wedge \dots \wedge \widehat{X_{i}} \wedge 
\dots \wedge \widehat{X_{j}}\wedge \dots \wedge X_{n}.
\end{array}$$  

We will now show that the Koszul resolution can be deformed.\\

\begin{theorem}~\label{koszul resolution}
Let $\goth a$ be a Lie algebra and let 
$\left (e_{1}, \dots, e_{n}\right )$ be a basis of $\goth a$. 
Denote by $C_{i,j}^{a}$ the structure constants of $\goth a$ with
respect to the basis $(e_{1},\dots,e_{n})$ so that we have 
$[e_{i},e_{j}]={\displaystyle \mathop
  \sum_{k=1}^{n}}C_{i,j}^{a}e_{a}$. 
Consider  $U_{h}({\goth a})$ a deformation of 
$U({\goth a})$ given under the form 
$$U_{h}({\goth a})\simeq 
\frac{T_{k[[h]]}\left ( 
{\displaystyle \mathop \oplus_{i=1}^{n}} k[[h]]
e_{i} \right )}{I}$$
where $I$ is the closure (in the $h$-adic topology) 
of the two sided ideal generated by the
relations 
$$e_{i}\otimes e_{j}-e_{j}\otimes e_{i}- 
g^{\otimes}_{i,j}(e_{1},\dots ,e_{n},h)$$
where $g_{i,j}$ satisfies the following : 
$$\begin{array}{l}
g_{i,j} \in k[X_{1}, \dots, X_{n}][[h]]\;\;{\rm and}\;\;
\partial _{X}(g_{i,j})\geq 1\\
g_{i,j}^{\otimes}(e_{1},\dots, e_{n})=
{\displaystyle \mathop \sum_{a=1}^{n}}C_{i,j}^{a}e_{a}
 \; {\rm mod}\; 
hT_{k[[h]]}\left ( {\displaystyle \mathop \oplus_{i=1}^{n}}k[[h]]e_{i}\right ).\\
\end{array} $$
 $k[[h]]$ is an $U_{h}({\goth a})$-module (called the trivial
 $U_{h}({\goth a})$-module) if we let the $e_{i}$'s act
 trivially. 
There exists a resolution of the trivial $U_{h}({\goth a})$-module
$k[[h]]$,  
$K_{h}=\left ( U_{h}({\goth a})\otimes _{k}\wedge^{\bullet} {\goth a}, 
{\partial}^{\bullet}_{h}\right )$,  such that $Gr K_{h}$ is the resolution of
the trivial $U({\goth a})[h]$-module $k[h]$. 
\end{theorem}

{\bf Remarks :}

1) Any quantized universal enveloping algebra,  $U_{h}({\goth a})$, has a
presentation as in the theorem because we might write it  
$U_{h}({\goth a})=\left ( U_{h}({\goth a})'\right )^{\vee}$.

2) The proof of the theorem gives an algorithm to construct the
resolution $K_{h}$. 

3) By theorem  ~\ref{koszul resolution}, we even get a filtered resolution of
the $FU_{h}({\goth g})$-module $k[[h]]$. \\

{\it Proof of the theorem ~\ref{koszul resolution}:}

We will prove by induction that on $q$ that one can construct 
$\partial ^{h}_{0}, \dots, \partial^{h}_{q}$ morphisms of
$U_{h}({\goth a})$-modules such that :\\

$\bullet \forall r\in [1,q], \partial_{r-1}^{h}\partial _{r}^{h}=0.$
 
$$ \begin{array}{l}
\bullet \partial _{r}^{h}(1 \otimes e_{p_{1}} \wedge \dots \wedge e_{p_{r}})=
{\displaystyle \mathop \sum_{i=1}^{r}}(-1)^{i-1}e_{p_{i}}\otimes 
e_{p_{1}}\wedge \dots \wedge \widehat{e_{p_{i}}}
\wedge \dots \wedge e_{p_{r}} \\
+ {\displaystyle \mathop \sum_{k<l}}
{\displaystyle \mathop \sum_{a}}(-1)^{k+l}
C_{p_{k},p_{l}}^{a}1 \otimes e_{a}\wedge e_{p_{1}}\wedge \dots \wedge 
\widehat{e_{p_{k}}}\wedge \dots \wedge \widehat{e_{p_{l}}}\wedge 
\dots \wedge e_{p_{r}}+
\alpha_{p_{1},\dots,p_{r}}
\end{array}$$
with 
$\alpha_{p_{1},\dots,p_{r}} \in hU_{h}({\goth a})
\otimes \wedge^{r-1}({\goth a}) $ so that $G\partial_{r}^{h}$ is the
qth differential of the Koszul complex of the trivial 
$U({\goth a})[h]$-module $k[h]$. From proposition 3.0.2, this implies that 
$\Ker \partial^{h}_{r-1}=\Im \partial^{h} _{r}$. 

We take  $\partial_{0}^{h} : U_{h}({\goth a}) \to k[[h]]$ to be the algebra
morphism determined by $\partial_{0}^{h}(e_{i})=0$. 

We take ${\partial}^{h}_{1} : U_{h}({\goth a})\otimes_{k} {\goth a} \to
U_{h}({\goth a})$ to be  the morphism of $U_{h}({\goth a})$-modules
determined by $\partial^{h}_{1}(u \otimes e_{i})=ue_{i}$.

One writes 
$$g_{i,j}(e_{1},\dots,e_{n},h)
={\displaystyle \mathop \sum_{a=1}^{n}}P_{i,j}^{a}e_{a}+
{\displaystyle \mathop \sum_{a=1}^{n}}C_{i,j}^{a}e_{a}$$
where the $P_{i,j}^{a}$'s are in $hU_{h}({\goth a})$. 

We look for a morphism of $U_{h}({\goth a})$-modules,
$\partial_{2}^{h}:  U_{h}({\goth a})\otimes_{k} {\wedge ^{2} \goth a} \to
U_{h}({\goth a})\otimes_{k}{\goth a}$   under the form 
$$\partial_{2}^{h}(1 \otimes e_{i}\wedge e_{j})= 
e_{i} \otimes e_{j}-e_{j} \otimes e_{i}-
{\displaystyle \sum_{a}}C_{i,j}^{a}1 \otimes e_{a}
-\alpha_{i,j}$$
 where $\alpha_{i,j}$ is in $hU_{h}({\goth a})\otimes {\goth a}$. 

One has 
$$\partial_{1}^{h}\left (  
e_{i} \otimes e_{j}-e_{j} \otimes e_{i}-
{\displaystyle \sum_{a}}C_{i,j}^{a}1 \otimes
e_{a}\right )=
{\displaystyle \mathop \sum_{a=1}^{n}}P_{i,j}^{a}e_{a}=
\partial _{1}^{h}\left ( 
{\displaystyle \mathop \sum_{a=1}^{n}}P_{i,j}^{a}\otimes e_{a}
\right ).$$
We might  take 
$$\partial_{2}^{h}(1 \otimes e_{i}\wedge e_{j})= 
e_{i} \otimes e_{j}-e_{j} \otimes e_{i}-
{\displaystyle \mathop \sum_{a=1}^{n}}C_{i,j}^{a}1 \otimes e_{a}
-{\displaystyle \mathop \sum_{a=1}^{n}}P_{i,j}^{a}1 \otimes e_{a}.$$
We have $\partial^{h}_{1}\circ \partial_{2}^{h}=0$

Let $q\geq 2$. 
Assume that $\partial_{0}^{h}, \partial_{1}^{h}, \dots
\partial_{q}^{h}$ are constructed 
and let  us construct 
$\partial_{q+1}^{h}$ as required. 

We look for 
$\partial _{q+1}^{h}(1 \otimes e_{p_{1}} \wedge \dots \wedge e_{p_{q+1}})$ 
under the form 
$$\begin{array}{l}
\partial _{q+1}^{h}(1 \otimes e_{p_{1}} \wedge \dots \wedge e_{p_{q+1}})=
{\displaystyle \mathop \sum_{i=1}^{q+1}}(-1)^{i-1}e_{p_{i}}\otimes 
e_{p_{1}}\wedge \dots \wedge \widehat{e_{p_{i}}}
\wedge \dots \wedge e_{p_{q+1}} \\
+ {\displaystyle \mathop \sum_{r<s}}
{\displaystyle \mathop \sum_{a}}(-1)^{r+s}
C_{p_{r},p_{s}}^{a}1 \otimes e_{a}\wedge e_{p_{1}}\wedge \dots \wedge 
\widehat{e_{p_{r}}}\wedge \dots \wedge \widehat{e_{p_{s}}}\wedge 
\dots \wedge e_{p_{q+1}}+
\alpha_{p_{1},\dots,p_{q+1}}
\end{array}$$
where $\alpha_{p_{1},\dots,p_{q+1}}$ is in 
$hU_{h}({\goth a})\otimes \wedge ^{q}({\goth a})$. The term 

$$\begin{array}{l}
 \partial_{q}^{h} \left ( 
{\displaystyle \mathop \sum_{i=1}^{q+1}}(-1)^{i-1}e_{p_{i}}\otimes 
e_{p_{1}}\wedge \dots \wedge \widehat{e_{p_{i}}}
\wedge \dots \wedge e_{p_{q+1}} \right )+\\
+ \partial_{q}^{h}\left ( {\displaystyle \mathop \sum_{r<s}}
{\displaystyle \mathop \sum_{a}}C_{p_{r},p_{s}}^{a}(-1)^{r+s}
1 \otimes e_{a}\wedge e_{p_{1}}\wedge \dots \wedge 
\widehat{e_{p_{r}}}\wedge \dots \wedge \widehat{e_{p_{s}}}\wedge 
\dots \wedge e_{p_{q+1}}\right )
\end{array}$$

equals 0 modulo $h$.  Hence it is in $hU_{h}({\goth a})\otimes \wedge ^{q-1}({\goth a})$.

As $\partial^{h}_{q-1}\partial^{h}_{q}=0$, it is  
in $hKer\partial_{q-1}^{h}=hIm\partial_{q}^{h}$. The existence
of $\alpha_{p_{1},\dots,p_{q+1}}$ follows. 
Hence we have constructed $\partial_{q+1}^{h}$ as required. 
The complex 
$K_{h}=\left (U_{h}({\goth a})\otimes \wedge^{\bullet} {\goth a} , 
\partial_{\bullet}^{h}\right )$ is a resolution of the trivial 
$U_{h}({\goth a})$-module 
$k[[h]]$.$\Box$    \\

\subsection{Quantum duality and deformation of the Koszul complex}

We may construct resolutions of the trivial $F_{h}[{\goth g}]$ and 
$F_{h}[{\goth g}]^{\vee}$-modules that respects the quantum duality. 

\begin{theorem}
Let ${\goth g}$ be a Lie bialgebra, $F_{h}[{\goth g}]$ a QFSHA such that 
${\displaystyle \frac {F_{h}[{\goth g}]}{hF_{h}[{\goth g}]}}$ is
  isomorphic to $F[{\goth g}]$ as a topological Poisson Hopf algebra 
and $F_{h}[{\goth g}]^{\vee}=U_{h}({\goth g}^{*})$ the
quantization of $U({\goth g}^{*})$ constructed from $F_{h}[{\goth g}]$
by the quantum duality principle. Let $\bar{x}_{1},\dots , \bar{x}_{n}$ be
elements of $F[{\goth g}]$ such that 
$F[{\goth g}]\simeq k[[\bar{x}_{1},\dots,\bar{x}_{n}]]$. Choose  
$x_{1}, \dots ,x_{n}$ elements of $F_{h}[{\goth g}]$ such that 
$x_{i}=\bar{x}_{i}\;{\rm mod}\; h$ and $\epsilon_{h}(x_{i})=0$. Then 
$U_{h}({\goth g}^{*})\simeq k[\check{x}_{1},\dots, \check{x}_{n}][[h]]$ 
with $\check{x}_{i}=h^{-1}x_{i}$. 
Let $( \epsilon_{1},\dots , \epsilon_{n})$ be a basis of ${\goth
  g}^{*}$ and let $C_{i,j}^{a}$ the structural constants of ${\goth
  g}^{*}$ with respect to this basis.  
We can  construct a resolution of the trivial 
$F_{h}[{\goth g}]$-module 
$K_{\bullet}^{h}=\left ( F_{h}[{\goth g}]\otimes \wedge
  {\goth g}^{*}, \partial _{q}^{h} \right )$   
of the form 
$$\begin{array}{l}
\partial_{q}^{h}
\left ( 1 \otimes \epsilon_{p_{1}}\wedge \dots \wedge \epsilon_{p_{q}} 
\right )= 
{\displaystyle \sum_{i=1}^{q}}(-1)^{i-1}x_{i}\otimes 
\epsilon_{p_{1}}\wedge \dots \wedge \widehat{\epsilon_{p_{i}}}
\wedge \dots \wedge \epsilon_{p_{q}}\\
+{\displaystyle \sum_{r<s}\sum_{a}}(-1)^{r+s}hC_{p_{r},p_{s}}^{a}
1 \otimes \epsilon_{a}\wedge \epsilon_{p_{1}}\wedge \dots \wedge 
\widehat{\epsilon_{p_{r}}}
\wedge \dots \wedge 
\widehat{\epsilon_{p_{s}}}\wedge \dots \wedge \epsilon_{p_{q}} \\
+{\displaystyle \sum_{t_{1},\dots ,t_{q-1}}} 
h \alpha_{p_{1}, \dots, p_{q}}^{t_{1}, \dots , t_{q-1}}\otimes 
\epsilon_{t_{1}}\wedge \dots \wedge \epsilon_{t_{q-1}} 
\end{array}$$ 
such that $\alpha_{p_{1}, \dots, p_{q}}^{t_{1}, \dots , t_{q-1}} \in 
I=\epsilon_{h}^{-1}( hk[[h]] )$. Set 
$$\check{\alpha}_{p_{1}, \dots, p_{q}}^{t_{1}, \dots , t_{q-1}}
(\check{x}_{1}, \dots , \check{x}_{n})=
\alpha_{p_{1}, \dots, p_{q}}^{t_{1}, \dots , t_{q-1}}
(x_{1}, \dots , x_{n}).$$
$\check{\alpha}_{p_{1}, \dots, p_{q}}^{t_{1}, \dots , t_{q-1}}$ is in 
$hk[\check{x}_{1}, \dots , \check{x}_{n}][[h]]$. 
Define the morphism of $U_{h}({\goth g}^{*})$-modules 
$\check{\partial}_{q}^{h} : 
U_{h}({\goth g}^{*})\otimes \wedge^{q}({\goth g}^{*}) \to 
U_{h}({\goth g}^{*})\otimes \wedge^{q-1}({\goth g}^{*})$ by 
$$\begin{array}{l}
\check{\partial}_{q}^{h}
\left ( 1 \otimes \epsilon_{p_{1}}\wedge \dots \wedge \epsilon_{p_{q}} 
\right )= 
{\displaystyle \sum_{i=1}^{n}}(-1)^{i-1}\check{x}_{i}\otimes 
\epsilon_{p_{1}}\wedge \dots \wedge \widehat{\epsilon_{p_{i}}}
\wedge \dots \wedge \epsilon_{p_{q}}\\
+{\displaystyle \sum_{r<s}\sum_{a}}(-1)^{r+s}C_{p_{r},p_{s}}^{a}
1 \otimes \epsilon_{a}\wedge \epsilon_{p_{1}}\wedge \dots \wedge 
\widehat{\epsilon_{p_{r}}}\wedge \dots \wedge \widehat{\epsilon_{p_{s}}}
\wedge \dots \wedge \epsilon_{p_{q}} \\
+ {\displaystyle \sum_{t_{1},\dots ,t_{q-1}}}
\check{\alpha}_{p_{1}, \dots, p_{q}}^{t_{1}, \dots , t_{q-1}}\otimes 
\epsilon_{t_{1}}\wedge \dots \wedge \epsilon_{t_{q-1}}. 
\end{array}$$ 
 Then $\check{K}_{h}^{\bullet}=\left ( U_{h}({\goth g}^{*})\otimes 
\wedge^{\bullet} {\goth g}^{*}, \check{\partial} _{q}^{h} \right )$ 
is a resolution of the trivial $U_{h}({\goth g}^{*})$-module.  
\end{theorem} 

{\it Proof of the theorem :}

One sets $x_{i}x_{j}-x_{j}x_{i}=
{\displaystyle \sum_{a=1}^{n}hC_{i,j}^{a}x_{a}+hu_{i,j}^{a}x_{a}}$. We
know that $u_{i,j}^{a}$ is in $I$.   
We take 
$$\begin{array}{l}
\partial_{0}^{h}=\epsilon_{h}\\
\partial_{1}^{h}(1 \otimes \epsilon _{i})=x_{i}\\
\end{array}$$
We set 
$$\partial_{2}^{h}(1 \otimes e_{i}\wedge e_{j})=x_{i}\otimes\epsilon_{j} 
- x_{j}\otimes \epsilon_{i}-
{\displaystyle \sum_{a}}hC_{i,j}^{a}\otimes \epsilon_{a}
-h{\displaystyle \sum_{a}}u_{i,j}^{a}\otimes \epsilon_{a}.$$
We have $\partial_{1}^{h} \circ \partial_{2}^{h}=0$ and we may choose 
$\alpha_{i,j}^{a}=u_{i,j}^{a}$.

Assume that 
$\partial^{h}_{0}, \partial^{h}_{1}, \dots , \partial^{h}_{q}$ have
been constructed such that 
\begin{itemize}
\item $\forall r \in [1,q] \;\;
\partial ^{h}_{r-1}\partial^{h}_{r}=0$
\item $\forall r \in [1,q] \;\;Im\partial^{h}_{r}=Ker
  \partial^{h}_{r-1}$ and satisfying the required relation.
\item $\alpha_{p_{1},p_{2}, \dots , p_{r}}^{q_{1},\dots ,q_{r-1}} \in I .$
\end{itemize}
and let us show that we can construct $\partial_{q+1}^{h}$ satisfying
these three conditions.

The computation
below is in [Kn] p 173.  

$$\begin{array}{l}
 \partial_{q}^{h} \left ( 
{\displaystyle \mathop \sum_{i=1}^{q+1}}(-1)^{i-1}x_{p_{i}}\otimes 
e_{p_{1}}\wedge \dots \wedge \widehat{\epsilon_{p_{i}}}
\wedge \dots \wedge \epsilon_{p_{q+1}} \right )+\\
+ \partial_{q}^{h}\left ( {\displaystyle \mathop \sum_{r<s}}
{\displaystyle \mathop \sum_{a}}hC_{p_{r},p_{s}}^{a}(-1)^{r+s}
1 \otimes \epsilon_{a}\wedge \epsilon_{p_{1}}\wedge \dots \wedge 
\widehat{\epsilon_{p_{r}}}\wedge \dots \wedge \widehat{\epsilon_{p_{s}}}\wedge 
\dots \wedge \epsilon_{p_{q+1}}\right )\\
= {\displaystyle \mathop \sum_{j<i}(-1)^{i+j}}
(x_{p_{i}}x_{p_{j}}-x_{p_{j}}x_{p_{i}})\otimes
\epsilon_{1}\wedge \dots \wedge 
\widehat{\epsilon_{p_{j}}}\wedge \dots \wedge \widehat {\epsilon_{p_{i}}}\wedge \dots
\wedge \epsilon_{p_{q+1}} \\
{\displaystyle \mathop \sum_{i}\sum_{r<s, r,s \neq i}}
{\displaystyle \mathop \sum_{a}}
(-1)^{r+s+\delta +i+1}hC_{p_{r},p_{s}}^{a}x_{p_{i}}
\otimes \epsilon_{a}\wedge \epsilon_{p_{1}}\wedge 
({\rm omit}\;p_{r},p_{s},p_{i})\wedge \epsilon_{p_{q+1}}\\
{\displaystyle \mathop \sum_{p_{r}<p_{s}}}
{\displaystyle \mathop \sum_{a}}
(-1)^{r+s}
hC_{p_{r},p_{s}}^{a}x_{a}\otimes \epsilon_{p_{1}}\wedge 
({\rm \; omit \; p_{r},p_{s}})\wedge
\epsilon_{p_{n+1}}\\
+ {\displaystyle \mathop \sum_{r<s}\sum_{a}\sum_{p_{i} \neq p_{r},p_{s}}}
(-1)^{r+s+i+\delta}
hC_{p_{r},p_{s}}^{a}x_{p_{i}}\otimes \epsilon_{a}\wedge \epsilon_{p_{1}}\wedge 
({\rm omit}\;p_{r},p_{s},p_{i})\wedge \epsilon_{p_{q+1}}\\
{\displaystyle \sum_{r<s}}(-1)^{r+s}
{\displaystyle \mathop \sum_{k<l;k,l\neq r,s}}(-1)^{k+l+\sigma}
{\displaystyle \mathop \sum_{a,b}}
h^{2}C_{p_{r},p_{s}}^{a}C_{p_{k},p_{l}}^{b}1 \otimes
\epsilon_{a}\wedge \epsilon_{b}\wedge \epsilon_{1}\wedge 
({\rm omit }\;p_{k},p_{l},p_{r},p_{s})\wedge \epsilon_{p_{q+1}}\\
+{\displaystyle \mathop \sum_{r<s}}(-1)^{r+s}
{\displaystyle \mathop \sum_{j\neq r,s}\sum_{a,b}}
(-1)^{j+\tau}
h^{2}C_{p_{r},p_{s}}^{a}C_{a,p_{j}}^{b}\otimes \epsilon_{b}\wedge \epsilon_{p_{1}}\wedge 
({\rm omit}\;p_{j},p_{r},p_{s})\wedge \epsilon_{p_{q+1}}\\
+ {\displaystyle \mathop \sum_{i}}(-1)^{i-1}
x_{p_{i}}h\alpha_{p_{1},\dots, \widehat{p_{i}},\dots,
  p_{q+1}}
+{\displaystyle \mathop \sum_{r<s}\sum_{a}}(-1)^{r+s} hC_{p_{r},p_{s}}^{a}
h\alpha_{a,p_{1},\dots, \widehat{p_{r}},\dots, \widehat{p_{s}},\dots,
  p_{q+1}}
  
\end{array}$$
where 
$$\begin{array}{l}
\delta =1\; {\rm if}\; r<i<s \; {\rm and}\; \delta =0 \;
{\rm otherwise}\\
{\sigma =1}{\; \rm if\;exactly\;one\;of\;k\;and\;l\;
  is\;between\;r\;and\;s}\\
\tau =1 \;{\rm if}\;r<j<s\;\;{\rm and}\; \tau =0\;{\rm otherwise}
\end{array}$$
The fifth term and the sixth  term cancel. The second
term and fourth term cancel with each other so that we have 
$$\begin{array}{l}
 \partial_{q}^{h} \left ( 
{\displaystyle \mathop \sum_{i=1}^{q+1}}\epsilon_{p_{i}}\otimes 
\epsilon_{p_{1}}\wedge \dots \wedge \widehat{\epsilon_{p_{i}}}
\wedge \dots \wedge \epsilon_{p_{q+1}} \right )+\\
+ \partial_{q}^{h}\left ( {\displaystyle \mathop \sum_{k<l}}
{\displaystyle \mathop \sum_{a}}
hC_{p_{k},p_{l}}^{a}1
 \otimes \epsilon_{a}\wedge \epsilon_{p_{1}}\wedge \dots \wedge 
\widehat{\epsilon_{p_{k}}}\wedge \dots \wedge \widehat{\epsilon_{p_{l}}}\wedge 
\dots \wedge \epsilon_{p_{q+1}}\right )\\
= {\displaystyle \mathop \sum_{j<i}(-1)^{i+j}}
\left (x_{p_{i}}x_{p_{j}}-x_{p_{j}}x_{p_{i}}
- {\displaystyle \mathop \sum_{a}}
hC_{p_{i},p_{j}}^{a}x_{a}\right )\otimes
\epsilon_{1}\wedge \dots \wedge 
\widehat{\epsilon_{p_{j}}}\wedge \dots \wedge \widehat {\epsilon_{p_{i}}}\wedge \dots
\wedge \epsilon_{p_{q+1}}\\
+ {\displaystyle \mathop \sum_{i}}(-1)^{i-1}
hx_{p_{i}}\alpha_{p_{1},\dots, \widehat{p_{i}},\dots,
  p_{q+1}}
+{\displaystyle \mathop \sum_{r<s}}(-1)^{r+s}h^{2}C_{p_{r},p_{s}}^{a}
\alpha_{a, p_{1},\dots, \widehat{p_{r}},\dots, \widehat{p_{l}},\dots,
  p_{q+1}}.\\
\end{array}$$
As $\partial^{h}_{q-1}\partial^{h}_{q}=0$, the term 
$$\begin{array}{l}
 \partial_{q}^{h} \left ( 
{\displaystyle \mathop \sum_{i=1}^{q+1}}(-1)^{i-1}e_{p_{i}}\otimes 
e_{p_{1}}\wedge \dots \wedge \widehat{e_{p_{i}}}
\wedge \dots \wedge e_{p_{q+1}} \right )+\\
+ \partial_{q}^{h}\left ( {\displaystyle \mathop \sum_{k<l}}
{\displaystyle \mathop \sum_{a}}(-1)^{k+l}
C_{p_{k},p_{l}}^{a}1 \otimes e_{a}\wedge e_{p_{1}}\wedge \dots \wedge 
\widehat{e_{p_{k}}}\wedge \dots \wedge \widehat{e_{p_{l}}}\wedge 
\dots \wedge e_{p_{q+1}}\right )
\end{array}$$
is in $hKer\partial_{q-1}^{h}=hIm\partial_{q}^{h}$.  We can choose 
$\alpha_{p_{1}, \dots ,p_{q+1}}^{t_{1},\dots ,t_{q}}$ in 
$F_{h}[{\goth g}]$ so that the expression above equals 
$-\partial_{q}^{h}
\left ( h \alpha_{p_{1}, \dots ,p_{q+1}}^{t_{1},\dots ,t_{q}} \right )$. 

Let us now prove that $\alpha_{p_{1}, \dots ,p_{q+1}}^{t_{1},\dots
  ,t_{q}}$ is in $I$. 
It is easy to see that \linebreak  
$-\partial_{q}^{h}\left (h \alpha_{p_{1}, \dots ,p_{q+1}}^{t_{1},\dots ,t_{q}} 
\otimes \epsilon_{t_{1}}\wedge \dots \wedge \epsilon_{t_{q}} \right )$
is element of $I^{3}\otimes \wedge^{q}{\goth g}^{*}$. 
Note that $\partial _{q}^{h}$ sends 
$I^{r}\otimes \wedge^{q}{\goth g}^{*}$ to 
$I^{r+1}\otimes \wedge^{q}{\goth g}^{*}$.
Let us write  
$$\alpha_{p_{1}, \dots ,p_{q+1}}^{t_{1},\dots ,t_{q}}= 
{\displaystyle \sum_{i_{1},\dots ,i_{n}}}
(\alpha_{p_{1}, \dots ,p_{q+1}}^{t_{1},\dots ,t_{q}})_{i_{1},\dots,i_{n}}
x_{1}^{i_{1}}\dots x_{n}^{i_{n}}$$
with $(\alpha_{p_{1}, \dots ,p_{q+1}}^{t_{1},\dots
  ,t_{q}})_{i_{1},\dots,i_{n}}$ in $k[[h]]$.
From the remarks we have just made, we see that 
$\partial_{q}^{h}\left ( h  {\displaystyle \sum_{t_{1},\dots ,t_{q}}}
(\alpha_{p_{1}, \dots ,p_{q+1}}^{t_{1},\dots ,t_{q}})_{0,\dots ,0}
\epsilon_{t_{1}}\wedge \dots \wedge \epsilon_{t_{q}}\right )$ is
in $I^{3}\otimes \wedge^{q}{\goth g}^{*}$. Hence 
$(\alpha_{p_{1}, \dots ,p_{q+1}}^{t_{1},\dots ,t_{q}})_{0,\dots ,0}$
is in $hk[[h]]$.

As ${\rm Im} G \partial_{q+1}^{h}={\rm Ker} G \partial_{q}^{h}$, one has 
${\rm Im}\partial_{q+1}^{h}={\rm Ker }\partial_{q}^{h}$. 

Set 
$$\check{\alpha}_{p_{1}, \dots , p_{q}}^{t_{1},\dots ,t_{q-1}}
(\check{x}_{1}, \dots , \check{x}_{n})=
\alpha_{p_{1}, \dots , p_{q}}^{t_{1},\dots ,t_{q-1}}
(x_{1}, \dots , x_{n}).$$
$$\begin{array}{l}
\check{\partial}_{0}=\epsilon \\
\check{\partial}_{1}(1 \otimes \epsilon_{i})=\check{x}_{i}\\
\check{\partial}_{2}(1 \otimes \epsilon_{i}\wedge \epsilon_{j})=
\check{x}_{i}\otimes \epsilon _{j}- \check{x}_{j}\otimes \epsilon_{j} 
-{\displaystyle \sum_{a}C_{i,j}^{a}}\otimes \epsilon_{a}-
{\displaystyle \sum_{a}\check{u}_{i,j}^{a}}\otimes \epsilon_{a}\\
\check{\partial}_{q+1}^{h}
\left ( 1 \otimes \epsilon_{p_{1}}\wedge \dots \wedge \epsilon_{p_{q+1}} 
\right )= 
{\displaystyle \sum_{i=1}^{q+1}}(-1)^{i-1}\check{x}_{i}\otimes 
\epsilon_{p_{1}}\wedge \dots \wedge \hat{\epsilon}_{p_{i}}\wedge \dots \wedge
\epsilon_{p_{q+1}}\\
+{\displaystyle \sum_{r<s}\sum_{a}}(-1)^{r+s}C_{p_{r},p_{s}}^{a}
1 \otimes \epsilon_{a}\wedge \epsilon_{p_{1}}\wedge \dots \wedge 
\hat{\epsilon}_{p_{r}}\wedge \dots \wedge \hat{\epsilon}_{p_{s}}
\wedge \dots \wedge \epsilon_{p_{q+1}} \\
+ {\displaystyle \sum_{t_{1}, \dots ,t_{q-1}}}
\check{\alpha}_{p_{1}, \dots, p_{q+1}}^{t_{1}, \dots , t_{q}}\otimes 
\epsilon_{t_{1}}\wedge \dots \wedge \epsilon_{t_{q}}. 
\end{array}$$
If $P$ is in $F_{h}$, one has 
$$\partial_{q}(P \otimes \epsilon_{p_{1}}\wedge \dots\wedge  \epsilon_{p_{q}})= h
\check{\partial }(\check{P} \otimes \epsilon_{p_{1}}\wedge \dots \wedge  \epsilon_{p_{q}}).$$
The relation $\check{\partial}_{q}\check{\partial}_{q+1}=0$ is
obtained by multiplying the relation  
$\partial_{q}^{h}\partial_{q+1}^{h}=0$
by $h^{-2}$. As $G \check{\partial}_{q}^{h}$ is the differential of the
Koszul complex of the trivial $U({\goth g}^{*})[h]$-module, the complex 
$\check{K}_{h}^{\bullet}=\left ( U_{h}({\goth g}^{*})\otimes 
\wedge^{\bullet} {\goth g}^{*}, \check{\partial} _{n}^{h} \right )$ 
is a resolution of the trivial $U_{h}({\goth g}^{*})$-module. $\Box $

\subsection{A link between $\theta_{F_{h}}$ and $\theta _{F_{h}^{\vee}}$}$\;$\\

\begin{theorem}~\label{link}
One has $\theta_{F_{h}}=h\theta_{F_{h}^{\vee}}$
\end{theorem}

{\it Proof of the theorem :}

We keep the notation of the previous
proposition and we will use the proof of the theorem ~\ref{qtrad}.  

The complex 
$\left ( \wedge^{\bullet}{\goth g}^{*}\otimes F_{h}, 
^{t}\partial_{n}^{h}\right )$ computes the $k[[h]]$-modules 
$\Ext^{i}_{F_{h}}\left (k[[h]], F_{h}\right )$. 
The cohomology class $cl  \left (1 \otimes \epsilon_{1}^{*}\wedge
  \dots \wedge \epsilon_{n}^{*} \right )$ is a  basis of  
$\underline{\Ext}_{F[{\goth g}][h]}^{n}\left (k[h], F[{\goth g}][h]\right )
\simeq G\Ext^{n}_{F_{h}}\left (k[[h]], F_{h}\right ) $. 
Hence there exists 
$\sigma =1+h\sigma_{1}+\dots \in \Ker ^{t}\partial _{n}^{h}$
such that  
$\left [cl \left ( \sigma \otimes \epsilon_{1}^{*}\wedge \dots 
\wedge \epsilon_{n}^{*}\right ) \right ]$ is a basis of 
$G\Ext^{n}_{F_{h}}\left (k[[h]], F_{h}\right )$. As the filtration on 
$\Ext^{n}_{F_{h}}\left (k[[h]], F_{h}\right )$ is Hausdorff, 
the cohomology class 
$cl\left (\sigma \otimes \epsilon_{1}^{*}\wedge \dots 
\wedge \epsilon_{n}^{*} \right )$ is a basis of 
$\Ext^{n}_{F_{h}}\left (k[[h]], F_{h}\right )$.

Define $\check{\sigma}$ by 
$$\check{\sigma}(\check{x}_{1}, \dots , \check{x}_{n})= 
\sigma (x_{1}, \dots , x_{n} ).$$
One has $^{t}\partial_{n}=h^{t}\check{\partial}_{n}$ and it is easy to
check that $\check{\sigma} \otimes \epsilon^{*}_{1}\wedge \dots \wedge 
\epsilon^{*}_{n}$ is in $\Ker ^{t}\check{\partial}^{h}_{n-1}$. 
If we had 
$$\check{\sigma}\otimes 
\epsilon _{1}^{*}\wedge \dots \wedge \epsilon_{n}^{*} =
^{t}\check{\partial}_{n-1}^{h}\left ( 
{\displaystyle \sum_{i=1}^{n}}\check{\sigma}_{i}\otimes 
\epsilon _{1}^{*} \wedge \dots \wedge 
\widehat{\epsilon_{i}^{*}} 
\wedge \dots \wedge  \epsilon_{n}^{*}
\right ),$$ 
then, reducing modulo $h$, we would get 
$$\overline{\check{\sigma}}\otimes 
\epsilon _{1}^{*}\wedge \dots \wedge \epsilon_{n}^{*} =
\overline{^{t}\check{\partial}_{n-1}^{h}}\left ( 
{\displaystyle \sum_{i=1}^{n}}\overline{\check{\sigma}_{i}}\otimes 
\epsilon _{1}^{*} \wedge \dots \wedge 
\widehat{\epsilon_{i}^{*}} 
\wedge \dots \wedge  \epsilon_{n}^{*}
\right ).$$
This would implies that 
$cl\left (1 \otimes \epsilon_{1}^{*}\wedge \dots \wedge \epsilon_{n}^{*}
 \right )$ is $0$ in 
$\Ext_{U({\goth g}^{*})}^{n}\left (k, U({\goth g}^{*})\right )$, 
which is impossible because 
$cl\left (1 \otimes \epsilon_{1}^{*}\wedge
  \dots \wedge \epsilon_{n}^{*} \right )$ is a  basis of  
$\Ext_{U({\goth g}^{*})}^{n}\left (k, U({\goth g}^{*})\right )$.
Thus 
$cl \left (\check{\sigma}\otimes 
\epsilon _{1}^{*}\wedge \dots \wedge \epsilon_{n}^{*} 
 \right )$ is a non zero element of  
$Ext^{dim {\goth g}^{*}}_{U_{h}({\goth g}^{*})}
\left ( k[[h]], U_{h}({\goth g}^{*})\right )$.  
For all $i$ in $[1,n]$, one has the relation 
$$ \sigma x_{i}\otimes 
\epsilon _{1}^{*}\wedge \dots \wedge \epsilon_{n}^{*} =
\theta _{F_{h}}(x_{i}) \sigma \otimes 
\epsilon _{1}^{*}\wedge \dots \wedge \epsilon_{n}^{*}+
^{t} \partial _{n}^{h}(\mu ) $$
Let us write 
$$\mu = {\displaystyle \sum_{i}}\mu_{i} \otimes 
\epsilon _{1}^{*} \wedge \dots \wedge
\widehat{\epsilon_{i}^{*}}
\wedge \dots \wedge  \epsilon_{n}^{*}$$
with $\mu_{i} \in F_{h}[{\goth g}]$. 
We set  $\check{\mu}_{i}(\check{x}_{1},\dots , \check{x}_{n} )=
\mu_{i}(x_{1}, \dots , x_{n})$ and 
$$\check{\mu} = {\displaystyle \sum_{i}}\check{\mu_{i}} \otimes 
\epsilon _{1}^{*} \wedge \dots \wedge
\widehat{\epsilon_{i}^{*}}
\wedge \dots \wedge  \epsilon_{n}^{*}.$$
Then we have 
$$ h\check{\sigma} \check{x}_{i}\otimes 
\epsilon _{1}^{*}\wedge \dots \wedge \epsilon_{n}^{*} =
\theta _{F_{h}}(x_{i}) \check{\sigma } \otimes 
\epsilon _{1}^{*}\wedge \dots \wedge \epsilon_{n}^{*}+
h^{t} \check{\partial} _{n}^{h} (\check{\mu} ).  $$
This finishes the proof of the theorem ~\ref{link}. $\Box$

\section{Study of on example}
We will now study explicitely an example suggested by B. Enriquez.
Chloup ([Chl]) introduced  the triangular Lie bialgebra 
$\left ({\goth g}= kX_{1}\oplus kX_{2} \oplus kX_{3} \oplus kX_{4} 
\oplus kX_{5} \right .$, 
\linebreak $\left . r=4 (X_{2}\wedge X_{3}) \right )$ where 
the non zero brackets are given by 
$$[X_{1},X_{2}]=X_{3},\;\; [X_{1},X_{3}]=X_{4},\;\; [X_{1}, X_{4}]=X_{5}$$
and the cobracket $\delta_{\goth g}$ is the following :
$$\forall X \in {\goth g},\;\; \delta (X)=X\cdot 4(X_{2} \wedge X_{3}).$$
 The dual Lie bialgebra of ${\goth g}$ will be denoted 
$\left ({\goth a}=
ke_{1}\oplus ke_{2} \oplus ke_{3} \oplus ke_{4} \oplus ke_{5}, \delta
\right )$. The only non zero Lie bracket of ${\goth a}$ is 
$[e_{2},e_{4}]=2e_{1}$ and its cobracket 
$\delta$ is non zero on the basis vectors $e_{3}$, $e_{4}$, $e_{5}$ :
$$\delta (e_{3})=e_{1}\otimes e_{2}-e_{2}\otimes e_{1}=
2 e_{1}\wedge e_{2},\;\;
\delta (e_{4})=2e_{1}\wedge e_{3},\;\;
\delta (e_{5})=2e_{1}\wedge e_{4}.$$
The invertible element of   
$U({\goth g})[[h]] \widehat {\otimes }U({\goth g})[[h]] $, 
$R=exp\left (h(X_{2} \otimes X_{3} - X_{3} \otimes X_{2}) \right )$, 
satisfies the equations 
$$\begin{array}{l}
R^{12}(\Delta \otimes 1)(R)=R^{23}(1 \otimes \Delta )(R)\\
(\epsilon \otimes id )(R)=1= (id \otimes \epsilon)(R).
\end{array}  $$
Thus, we may twist the trivial deformation of 
$\left ( U({\goth g})[[h]], \mu _{0}, \Delta_{0}, \iota_{0},\epsilon_{0} ,S_{0} 
\right )$ by 
$R$ ([C-P] p. 130). The topological Hopf algebra obtained has the same 
multiplication, antipode , unit and counit  but its coproduct is
$\Delta ^{R}=R^{-1}\Delta_{0}R$. It is a quantization of  
 $\left ( {\goth g}, r \right )$. 
We will denote it by $U_{h}({\goth  g})$.
The Hopf algebra 
$U_{h}({\goth g})^{*}$ is a QFSHA and  
$\left ( U_{h}({\goth g})^{*}\right ) ^{\vee} $ is a quantization of 
$\left ( {\goth a}, \delta_{\goth a} \right ) $. We will compute it
  explicitely.

\begin{proposition}\label{example}
a) $\left ( U(\goth g)^{*}\right ) ^{\vee}$ is isomorphic as a
topological Hopf algebra to the topological k[[h]]-algebra 
$T_{k[[h]]}\left ( k[[h]]e_{1}\oplus k[[h]]e_{2}\oplus k[[h]]e_{3}
\oplus k[[h]]e_{4} \oplus k[[h]]e_{5} \right )/I$
where $I$ is the closure of the two-sided ideal generated by 
$$\begin{array}{l}
e_{2}\otimes e_{4}-e_{4}\otimes e_{2}-2e_{1}\\
e_{3}\otimes e_{5}-e_{5}\otimes e_{3}-
{\displaystyle \frac{2}{3}}h^{2}e_{1}\otimes e_{1}\otimes e_{1}\\
e_{4}\otimes e_{5}-e_{5}\otimes e_{4}-
{\displaystyle \frac{1}{6}}h^{3}
e_{1}\otimes e_{1}\otimes e_{1}\otimes e_{1}\\
e_{2}\otimes e_{5}-e_{5}\otimes e_{2}+he_{1}\otimes e_{1}\\
e_{3}\otimes e_{4}-e_{4}\otimes e_{3}+he_{1}\otimes e_{1}\\
e_{i}\otimes e_{j}-e_{j}\otimes e_{i}\;{\rm if}\; 
\{i,j \}\neq \{2,4\},\; \{3,5\}, \; \{4,5\},\; \{2,5\},\; \{3,4\}
\end{array}$$
with the coproduct $\Delta_{h}$, counit $\epsilon_{h}$ and antipode
$S$ defined as follows :
$$\begin{array}{l}
\Delta_{h}(e_{1})=e_{1}\otimes 1 + 1 \otimes e_{1}\\
\Delta_{h}(e_{2})=e_{2}\otimes 1 + 1 \otimes e_{2}\\
\Delta_{h}(e_{3})=e_{3}\otimes 1 +1 \otimes e_{3}-he_{2}\otimes
e_{1}\\
\Delta_{h}(e_{4})=e_{4}\otimes 1 + 1 \otimes e_{4}-he_{3}\otimes e_{1}
+{\displaystyle \frac {h^{2}}{2}}e_{2}\otimes e_{1}^{2}\\
\Delta_{h} (e_{5})=e_{5}\otimes 1 +1\otimes e_{5}-he_{4}\otimes e_{1}
+{\displaystyle \frac{h^{2}}{2}}e_{3}\otimes e_{1}^{2}
-{\displaystyle \frac{h^{3}}{6}}e_{2}\otimes e_{1}^{3}.\\
\forall i \in [1,5],\;\;\epsilon_{h}(e_{i})=0\\
\forall i \in [1,5],\;\; S(e_{i})=-e_{i}
\end{array}$$
b) $\left ( U(\goth g)^{*}\right ) ^{\vee}$ is not isomorphic to the
trivial deformation of $U({\goth a})$, $U({\goth a})[[h]]$, 
 as  algebra. 
\end{proposition}

{\it Proof of the proposition}

Let $\xi_{i}$ be the element of $U({\goth g})^{*}$ defined by 
$$<\xi_{i},X^{a_{1}}_{1}X_{2}^{a_{2}}X_{3}^{a_{3}}X_{4}^{a_{4}}X_{5}^{a_{5}}>=
\delta_{a_{1},0}\dots \delta_{a_{i},1}\dots \delta_{a_{5},0}.$$
 
The algebras $U({\goth g})^{*}$ and $k[[\xi_{1}, \dots, \xi_{n}]]$ are
isomorphic. The topological Hopf algebra 
$\left ( U_{h}({\goth g}) ^{*}, ^{t}\Delta_{0}^{R}=\cdot_{h},
^{t}\mu_{0}=\Delta_{h}, ^{t}\epsilon_{0}, ^{t}\iota _{0}=\epsilon_{h} , 
^{t}S_{0}\right )$ is
a QFSHA. 
Remark that  $U_{h}({\goth g}) ^{*}$ and 
$k[[\xi_{1}, \dots, \xi_{n} , h]]$ are isomorphic as $k[[h]]$-modules. 
The elements  $\xi_{1}, \dots, \xi_{n}$ generate topologically  the $k[[h]]$-
algebra $ U_{h}({\goth g}) ^{*}$ and satisfy 
$\epsilon_{h}(\xi_{i})=0$.
$$< \xi_{2} \otimes \xi_{4} -\xi_{4}\otimes \xi_{2},
\Delta ^{R} (X_{1}^{a_{1}}\dots X_{5}^{a_{5}}) > \neq 0  \Longleftrightarrow 
(a_{1},a_{2}, a_{3},a_{4},a_{5})=(1,0,0,0,0).$$
and $< \xi_{2} \otimes \xi_{4} -\xi_{4}\otimes \xi_{2},\Delta ^{R} (X_{1})>
=2h$.
Hence  $\xi_{2}\cdot _{h} \xi_{4}-\xi_{4}\cdot_{h}\xi_{2}=2h
\xi_{1}$. 
$$< \xi_{3} \otimes \xi_{5} -\xi_{5}\otimes \xi_{3},
\Delta ^{R}(X_{1}^{a_{1}}\dots X_{5}^{a_{5}}) > \neq 0  \Longleftrightarrow 
(a_{1},a_{2}, a_{3},a_{4},a_{5})=(3,0,0,0,0)$$
and $< \xi_{3} \otimes \xi_{5} -\xi_{5}\otimes \xi_{3},X_{1}^{3}>=4h$.
Hence  $\xi_{3}\cdot _{h} \xi_{5}-\xi_{5}\cdot_{h}\xi_{3}=
{\displaystyle \frac{2h^{2}}{3}}
\xi_{1}\cdot_{h} \xi_{1}\cdot _{h} \xi_{1}$. 
$$< \xi_{4} \otimes \xi_{5} -\xi_{5}\otimes \xi_{4},
\Delta ^{R}(X_{1}^{a_{1}}\dots X_{5}^{a_{5}}) > \neq 0  \Longleftrightarrow 
(a_{1},a_{2}, a_{3},a_{4},a_{5})=(4,0,0,0,0).$$
and $< \xi_{4} \otimes \xi_{5} -\xi_{5}\otimes \xi_{4},
\Delta (X_{1}^{4})>=-4h$.
Hence  $\xi_{4}\cdot _{h} \xi_{5}-\xi_{5}\cdot_{h}\xi_{4}=
{\displaystyle \frac{-h^{3}}{6}}
\xi_{1}\cdot_{h} \xi_{1}\cdot_{h} \xi_{1}\cdot_{h} \xi_{1}$. 
$$< \xi_{2} \otimes \xi_{5} -\xi_{5}\otimes \xi_{2},
\Delta ^{R} (X_{1}^{a_{1}}\dots X_{5}^{a_{5}}) > \neq 0  \Longleftrightarrow 
(a_{1},a_{2}, a_{3},a_{4},a_{5})=(2,0,0,0,0).$$
and $< \xi_{2} \otimes \xi_{5} -\xi_{5}\otimes \xi_{2},
\Delta (X_{1}^{2})>=-2h$.
Hence  $\xi_{2}\cdot _{h} \xi_{5}-\xi_{5}\cdot_{h}\xi_{2}= -h 
\xi_{1}\cdot_{h} \xi_{1}$. 
$$< \xi_{3} \otimes \xi_{4} -\xi_{4}\otimes \xi_{3},
\Delta ^{R} (X_{1}^{a_{1}}\dots X_{5}^{a_{5}}) > \neq 0  \Longleftrightarrow 
(a_{1},a_{2}, a_{3},a_{4},a_{5})=(2,0,0,0,0).$$
and $< \xi_{3} \otimes \xi_{4} -\xi_{4}\otimes \xi_{3},
\Delta ^{R}(X_{1}^{2})>=-2h$.
Hence  $\xi_{3}\cdot _{h} \xi_{4}-\xi_{4}\cdot_{h}\xi_{3}=
-h  \xi_{1}\cdot_{h} \xi_{1}$. In the cases different from those
mentionned above, $\xi_{i}\cdot _{h}
\xi_{j}=\xi_{j}\cdot_{h}\xi_{i}$. 

Let us now compute the coproduct $\Delta_{h}$ of $U_{h}({\goth
  g})^{*}$. 
$$\begin{array}{l}
<\Delta_{h} (\xi_{3}), 
X_{1}^{a_{1}}X_{2}^{a_{2}}X_{3}^{a_{3}}X_{4}^{a_{4}}X_{5}^{a_{5}}
\otimes X_{1}^{b_{1}}X_{2}^{b_{2}}X_{3}^{b_{3}}X_{4}^{b_{4}}X_{5}^{b_{5}}>
\neq 0 \Longleftrightarrow \\
(a_{1}, a_{2}, a_{3}, a_{4}, a_{5}, b_{1}, b_{2}, b_{3}, b_{4}, b_{5})
=(0,0,1,0,0,0,0,0,0,0)\;{\rm or }\; (0,0,0,0,0,0,0,1,0,0) \\
{\rm or}\;\; (0,1,0,0,0, 1,0,0,0,0)
\end{array}$$
and $<\Delta_{h} (\xi_{3}), X_{2}X_{1}>= -1.$
Hence 
$$ \Delta_{h} (\xi_{3})= \xi_{3}\otimes 1 + 1 \otimes
\xi_{3}-\xi_{2}\otimes \xi_{1}.$$
$$\begin{array}{l}
<\Delta_{h} (\xi_{4}), 
X_{1}^{a_{1}}X_{2}^{a_{2}}X_{3}^{a_{3}}X_{4}^{a_{4}}X_{5}^{a_{5}}
\otimes X_{1}^{b_{1}}X_{2}^{b_{2}}X_{3}^{b_{3}}X_{4}^{b_{4}}X_{5}^{b_{5}}>
\neq 0 \Longleftrightarrow \\
(a_{1}, a_{2}, a_{3}, a_{4}, a_{5}, b_{1}, b_{2}, b_{3}, b_{4}, b_{5})
=(0,0,0,1,0,0,0,0,0,0)\;{\rm or }\; (0,0,0,0,0,0,0,0,1,0) \\
{\rm or}\;\; (0,0,1,0,0, 1,0,0,0,0)\;{\rm or}\; 
(0,1,0,0,0,2,0,0,0,0).
\end{array}$$
Moreover  
$$<\Delta _{h} (\xi_{4}), X_{3} \otimes X_{1}>= -1 \;\;{\rm and}\;\;
<\Delta_{h} (\xi_{4}), X_{2}\otimes X_{1}^{2}>= 1.$$
Hence 
$$ \Delta _{h}(\xi_{4})= \xi_{4}\otimes 1 + 1 \otimes
\xi_{4}-\xi_{3}\otimes \xi_{1}+
{\displaystyle \frac{1}{2}}\xi_{2} \otimes \xi_{1}\cdot_{h}\xi_{1}.$$
$$\begin{array}{l}
<\Delta_{h} (\xi_{5}), 
X_{1}^{a_{1}}X_{2}^{a_{2}}X_{3}^{a_{3}}X_{4}^{a_{4}}X_{5}^{a_{5}}
\otimes X_{1}^{b_{1}}X_{2}^{b_{2}}X_{3}^{b_{3}}X_{4}^{b_{4}}X_{5}^{b_{5}}>
\neq 0 \Longleftrightarrow \\
(a_{1}, a_{2}, a_{3}, a_{4}, a_{5}, b_{1}, b_{2}, b_{3}, b_{4}, b_{5})
=(0,0,0,0,1,0,0,0,0,0)\;{\rm or }\; (0,0,0,0,0,0,0,0,0,1) \\
{\rm or}\;\; (0,0,0,1,0, 1,0,0,0,0)\;{\rm or}\; 
(0,0,1,0,0,2,0,0,0,0)\;{\rm or}\; (0,1,0,0,0,3,0,0,0,0).
\end{array}$$
Moreover 
$$<\Delta_{h} (\xi_{5}), X_{4} \otimes X_{1}>= -1,\;\; 
<\Delta_{h} (\xi_{4}), X_{3}\otimes X_{1}^{2}>= 1,\;\; 
<\Delta_{h} (\xi_{4}), X_{2}\otimes X_{1}^{3}>= -1.$$
Hence 
$$ \Delta_{h} (\xi_{5})= \xi_{5}\otimes 1 + 1 \otimes
\xi_{5}-\xi_{4}\otimes \xi_{1}+
{\displaystyle \frac{1}{2}}\xi_{3} \otimes \xi_{1}\cdot_{h}\xi_{1}
-{\displaystyle \frac{1}{6}}\xi_{2}\otimes
\xi_{1}\cdot_{h}\xi_{1}\cdot_{h} \xi_{1}.$$

We set $\check{\xi_{i}}=h^{-1}\xi_{i}$ and 
$e_{i}=\check{\xi_{i}} \;{\rm mod}\; h \left ( U({\goth g})^{*}\right )^{\vee}$.
Let $\chi : \left ( U({\goth g})^{*}\right )^{\vee} \to U({\goth a})[[h]]$ 
be the isomorphism of topologicall $k[[h]]$-modules defined by 
$$\chi \left ( {\displaystyle \mathop \sum_{r \in \N}}
P_{r}(\check{\xi} _{1}, \dots ,\check{\xi}_{n})h^{r}  \right )= 
{\displaystyle \mathop \sum_{r \in \N}}
P_{r}(e_{1}, \dots ,e_{n})h^{r}.$$
\\
From what we have reviewed in the first paragraph of this section, the
first part of this theorem is proved. 

If $u$ and $v$ are in $U({\goth a})$, one sets 
$$u \cdot_{h} v =uv+{\displaystyle \mathop \sum_{r=1}^{\infty}}h^{r}
\mu_{r}(u,v).$$
one has 
$$\mu_{1}(e_{3},e_{4})=0,\;\; \mu_{1}(e_{4},e_{3})=e_{1}^{2},\;\;
\mu_{1}(e_{2},e_{5})=0,\;\; \mu_{1}(e_{5},e_{2})=e_{1}^{2}.$$
Let us show now that $\mu_{1}$ is a coboundary in the Hochschild
cohomology. The Hochschild cohomology 
$HH^{*}(U({\goth a}), U({\goth a}))$ 
is computed by the complex 
$\left ( Hom \left ( U({\goth a})^{\otimes}, U({\goth a})\right ), b \right )$
where : if 
$f \in Hom \left ( U({\goth a})^{\otimes n+1}, U({\goth a})\right )$,
then 
$$b (f) (a_{0}, \dots, a_{n})=a_{0}f(a_{1},\dots, a_{n})+ 
{\displaystyle \mathop \sum_{i=1}^{n}}(-1)^{i}f(a_{0}, \dots, a_{i-1}a_{i},
\dots a_{n})+ f(a_{0}, \dots,a_{n-1})a_{n}(-1)^{n}.$$
  The Lie algebra cohomology of ${\goth a}$ with coefficients in
  $U({\goth a})^{ad}$ (with the adjoint action),
$H^{*}\left ( {\goth a}, U({\goth a})^{ad} \right )$, is computed by
the Chevalley-Eilenberg complex 
$\left ( Hom \left ( \bigwedge {\goth a}, U({\goth a})\right ), d \right )$
where : if 
$f \in Hom \left ( \bigwedge^{ n+1}{\goth a}, U({\goth a})\right )$
$$\begin{array}{rcl}
d(f)(z_{1}, \dots,z_{n+1}) &=&
{\displaystyle \mathop \sum_{i=1}^{n+1}}(-1)^{i-1}z_{i} \cdot
f(z_{1}, \dots, z_{i-1},z_{i+1},\dots, z_{n+1})\\
&+&{\displaystyle \mathop \sum_{i<j}}(-1)^{i+j}
f([z_{i},z_{j}], \dots, z_{i-1}, z_{i+1}, z_{j-1},z_{j+1}, z_{n+1}).
\end{array}$$
The map  ([L] lemma 3.3.3)
$\Psi ^{*} : 
\left ( Hom \left ( U({\goth a})^{\otimes}, U({\goth a})\right ), b
\right )\to 
\left ( Hom \left ( \bigwedge {\goth a}, U({\goth a})^{ad}\right ), d \right )$
defined by antisymmetrization
$$
\Psi^{*}(f)(z_{1}, \dots, z_{n})  =  
f\left ( {\displaystyle \mathop \sum_{\sigma \in {S_{n}}}} \epsilon(\sigma ) 
z_{\sigma (1)} \otimes \dots \otimes z_{\sigma (n)} \right )$$
is a morphism of complexes. 
One checks easily that  
$$\Psi^{*} (\mu_{1})=d\left ( 
{\displaystyle -\frac {1}{2}e_{1}e_{2}\otimes e_{3}^{*}
-\frac {1}{2}e_{1}e_{4}\otimes e_{5}^{*}}
\right ).$$ 
There exists $\alpha \in Hom \left ( U({\goth a}), U({\goth a})\right )$
such that $\mu_{1}=b(\alpha)$. The map $\alpha$ is determined by 
$$\begin{array}{l}
\alpha_{\mid {\goth a}}= {\displaystyle -\frac {1}{2}e_{1}e_{2}\otimes e_{3}^{*}
-\frac {1}{2}e_{1}e_{4}\otimes e_{5}^{*}}\\
\forall (u,v) \in U({\goth a}),\;\; 
\mu_{1}(u,v)=u\alpha(v)-\alpha(uv)+u\alpha (v)
\end{array}$$ 

We set $\beta_{h}=id -h \alpha$. Then one has  
$\beta_{h}^{-1}={\displaystyle \mathop
  \sum_{i=0}^{\infty}}h^{i}\alpha^{i}$. If $u$ and $v$ are elements of  
$U({\goth a})$, we put  
$$u\cdot_{h}^{'}v= 
\beta_{h}^{-1}\left ( \beta_{h}(u)\cdot_{h}\beta_{h}(v)\right ).$$
Let's compute $e_{i}\cdot_{h}'e_{j}- e_{j}\cdot_{h}'e_{i}$. 
If $i$ and $j$ are different from
$3$ and $5$, then  $e_{i}\cdot_{h}'e_{j}=e_{i}\cdot_{h}e_{j}$
$$\begin{array}{rcl}
e_{1}\cdot_{h}'e_{3}-e_{3}\cdot_{h}'e_{1} & = &
f_{h}^{-1}
\left [ e_{1}\cdot_{h}\left ( e_{3}+
{\displaystyle \frac{he_{1}e_{2}}{2}}\right )- 
\left ( e_{3}+{\displaystyle \frac{he_{1}e_{2}}{2}} \right )
\cdot_{h}e_{1}\right ]\\
&=& f_{h}^{-1}\left [ 
e_{1}\cdot_{h}{\displaystyle \frac {he_{1}e_{2}}{2}}
+ {\displaystyle \frac {he_{1}e_{2}}{2}}\cdot_{h}e_{1}\right ]\\
&=& 0
\end{array}$$
Similarly, the following relations hold
$$\begin{array}{l}
e_{1}\cdot_{h}'e_{5}=e_{5}\cdot_{h}'e_{1},\;\; 
e_{2}\cdot_{h}'e_{3} =e_{3}\cdot_{h}'e_{2},\;\;
e_{2}\cdot_{h}' e_{5}=e_{5}\cdot_{h}'e_{2},\;\;
e_{3}\cdot_{h}'e_{4}=e_{4}\cdot_{h}'e_{3},
\end{array}$$

Let us now compute $e_{3}\cdot_{h}'e_{5}-e_{5}\cdot_{h}'e_{3}$. Easy
computations lead to the following equalities :
one has
$$\begin{array}{l}
e_{1}e_{2}\cdot_{h}e_{5}-e_{5}\cdot_{h}e_{1}e_{2}=e_{1}^{3}\\
e_{3}\cdot_{h}e_{1}e_{4}-e_{1}e_{4}\cdot_{h}e_{3}=-e_{1}^{3}\\
e_{1}e_{2}\cdot_{h}e_{1}e_{4}-e_{1}e_{4}\cdot_{h}e_{1}e_{2}=2e_{1}^{3}
\end{array}$$
One deduces easily from this that 
$$e_{3}\cdot_{h}^{'}e_{5}-e_{5}\cdot^{'}_{h}e_{3}=
{\displaystyle \frac{1}{6}}h^{2}e_{1}^{3}.$$
Similarly, one has 
$$e_{4}\cdot_{h}' e_{5}-e_{5}\cdot_{h}'e_{4}=
{\displaystyle \frac{-h^{2}}{6}}e_{1}^{3}.$$
The topological algebras 
$\left [U({\goth a})[[h]], \cdot_{h} \right ]$ and 
$\left [ U({\goth a})[[h]], \cdot_{h}' \right ]$ are isomorphic, hence
their centers are isomorphic. 
Let us compute the center of 
$\left [ U({\goth a})[[h]], \cdot_{h}' \right ]$.
Let $z$ be an element of the center 
$Z\left [ U({\goth a})[[h]], \cdot_{h}'\right ]$. One writes $z$ under
the form 
${\displaystyle \mathop \sum_{n\geq 0}}P_{r}(e_{1},e_{2},e_{3},e_{4},e_{5})h^{r}$
(where the multiplications in $P_{r}(e_{1},e_{2},e_{3},e_{4},e_{5})$
are $\cdot_{h}'$). One has 
$$e_{2}\cdot_{h}'z -z \cdot_{h}' e_{2}=
{\displaystyle \mathop \sum _{r \in \N}}2h^{r}
\frac{\partial P_{r}}{\partial X_{4}}(e_{1},e_{2},e_{3},e_{4},e_{5}).$$
Hence the polynomials $P_{r}$ don't depend on $X_{4}$ and $z$ can be
written 
$z={\displaystyle \mathop \sum_{n\geq
    0}}P_{r}(e_{1},e_{2},e_{3},e_{5})h^{r}$. 
$$e_{3}\cdot_{h}'z -z \cdot_{h}' e_{3}=
{\displaystyle \mathop \sum _{r \in \N}}\frac{1}{6}h^{r+2}
\left ( X_{1}^{3}\frac{\partial P_{r}}{\partial X_{5}}\right )
(e_{1},e_{2},e_{3},e_{5}).$$
Hence the polynomials $P_{r}$ don't depend on $X_{5}$ and $z$ can be
written 
$z={\displaystyle \mathop \sum_{n\geq 0}}P_{r}(e_{1},e_{2},e_{3})h^{r}$. 

$$e_{4}\cdot_{h}'z -z \cdot_{h}' e_{4}=
{\displaystyle \mathop \sum _{r \in \N}}-2h^{r}
\frac{\partial P_{r}}{\partial
  X_{2}}(e_{1},e_{2},e_{3}).$$
Hence the polynomials $P_{r}$ don't depend on $X_{2}$ and $z$ can be
written 
$z={\displaystyle \mathop \sum_{n\geq 0}}P_{r}(e_{1},e_{3})h^{r}$. 
$$e_{5}\cdot_{h}'z -z \cdot_{h}' e_{5}=
{\displaystyle \mathop \sum _{r \in \N}}\frac{-1}{6}h^{r+2}
\left ( X_{1}^{3}\frac{\partial P_{r}}{\partial X_{3}}\right )
(e_{1},e_{3}).$$
Hence the polynomials $P_{r}$ don't depend on $X_{3}$ and $z$ can be
written 
$z={\displaystyle \mathop \sum_{n\geq 0}}P_{r}(e_{1})h^{r}$.
Hence 
$$Z\left [ U({\goth a})[[h]], \cdot_{h}' \right ]=
\{{\displaystyle \mathop \sum_{n\geq 0}}P_{r}(e_{1})h^{r}\mid
P_{r}\in k[X_{1}]\}.$$
But, the center of the trivial deformation of $U({\goth a })$ is 
$$Z\left [ U({\goth a})[[h]], \mu_{0} \right ]=
\{{\displaystyle \mathop \sum_{n\geq 0}}P_{r}(e_{1},e_{3},e_{5})
h^{r}\mid P_{r}\in k[X_{1},X_{3},X_{5}]\}.$$
The algebras $\left [ U({\goth a})[[h]], \cdot_{h}' \right ]$ and 
$\left [ U({\goth a})[[h]], \mu_{0} \right ]$ are not isomorphic as
their center are not isomorphic. $\Box $

%{\bf Remarques :}

%En fait je crois qu'il serait plus inter{\'e}ssant de montrer que 
%$\left ( U_{h}({\goth g})^{*}\right )^{\vee}$ n'est pas isomorphe en tant
%qu'alg{\`e}bre {\`a} la  compl{\'e}tion de l'alg{\`e}bre enveloppante d'une
%d{\'e}formation de $\goth a$. Mais c'est plus compliqu{\'e}. 

%En effet reprenons la big{\`e}bre duale  de la big{\`e}bre de 
%Chloup $\goth a$. On d{\'e}finit la 
%$\C[[h]]$-alg{\`e}bre de Lie ${\goth a}[[h]]$ comme suit : les seuls
%crochets non nuls sont  
%$$[e_{3},e_{5}]=he_{3},\;\;\;
%[e_{2},e_{4}]=e_{1}$$
%$\widehat{U_{\C[[h]]}\left ( {\goth a}[[h]]\right )}$ est une
%quantification de $U({\goth a })$ qui n'est pas isomorphe {\`a} 
%$U({\goth a })[[h]]$ (quantification triviale) en tant qu'alg{\`e}bre. %

%Dans l'exemple {\'e}tudi{\'e} ci-dessus, je pense que l'on peut
%d{\'e}montrer, par la m{\^e}me m{\'e}thode, que 
%$\left [ U ({\goth a})[[h]], \cdot_{h}\right ]\simeq  
%\left [ U ({\goth a})[[h]], \cdot_{h}'\right ]$ est isomorphe {\`a}  
%$\left [ U ({\goth a})[[h]], \cdot_{h}''\right ]$  avec
%$$\begin{array}{l}
%e_{1}\cdot_{h}''e_{5}=e_{5}\cdot_{h}''e_{1},\;\; 
%e_{2}\cdot_{h}''e_{3} =e_{3}\cdot_{h}''e_{2},\;\;
%e_{2}\cdot_{h}'' e_{5}=e_{5}\cdot_{h}''e_{2},\;\;
%e_{3}\cdot_{h}''e_{4}=e_{4}\cdot_{h}''e_{3},\\
%e_{4}\cdot _{h}''e_{5}=e_{5}\cdot_{h}''e_{4}\\
%e_{3}\cdot_{h}''e_{5} - e_{5}\cdot_{h}''e_{3}=
%{\displaystyle \frac{h^{2}}{6}}e_{1}^{3}.\\
%\end{array}$$

\begin{proposition}
We consider the quantized enveloping algebra of the proposition \ref{example} 
We write the relations defining the ideal $I$ as follows 
$$e_{i}\otimes e_{j}-e_{j}\otimes e_{i}-
{\displaystyle \mathop \sum_{a}}C_{i,j}^{a}e_{a}
-P_{i,j}.$$
As all the $P_{i,j}$'s are monomials in $e_{1}$'s, the notation 
${\displaystyle \mathop \frac{P_{i,j}}{e_{1}}}$ makes sense. The
complex 
$$0 \to U_{h}({\goth a})\otimes \wedge ^{5}{\goth a}
\buildrel {\partial_{5}^{h}}\over \rightarrow 
U_{h}({\goth a})\otimes \wedge ^{4}{\goth a}
\buildrel {\partial_{4}^{h}} \over \rightarrow  \dots 
\buildrel {\partial_{2}^{h}}\over \rightarrow 
U_{h}({\goth a})\otimes {\goth a}
\buildrel {\partial_{1}^{h}}\over \rightarrow
U_{h}({\goth a})\buildrel {\partial_{0}^{h}}\over \rightarrow
k[[h]]\to 0$$
where the morphisms of  $U_{h}({\goth a})$, $\partial_{h}^{i}$, are
described below is a resolution of the trivial $U_{h}({\goth a})$-module 
$k[[h]]$.We set 
$$\begin{array}{l}
\partial _{n}(1 \otimes e_{p_{1}} \wedge \dots \wedge e_{p_{n}})=
{\displaystyle \mathop \sum_{i=1}^{n}}(-1)^{i-1}e_{p_{i}}\otimes 
e_{p_{1}}\wedge \dots \wedge \widehat{e_{p_{i}}}
\wedge \dots \wedge e_{p_{n}} \\
+ {\displaystyle \mathop \sum_{k<l}}(-1)^{k+l}
{\displaystyle \mathop \sum_{a}}
C_{p_{k},p_{l}}^{a}1 \otimes e_{a}\wedge e_{p_{1}}\wedge \dots \wedge 
\widehat{e_{p_{k}}}\wedge \dots \wedge \widehat{e_{p_{l}}}\wedge 
\dots \wedge e_{p_{n}}.
\end{array}$$
Then
$$\begin{array}{l}
\partial_{0}^{h}=\epsilon_{h}\\
\partial_{1}^{h}(1 \otimes e_{i})=e_{i}\\
\partial_{2}^{h}(1 \otimes e_{i}\wedge e_{j})=
\partial_{2}(1 \otimes e_{i}\wedge e_{j})- 
{\displaystyle \mathop \frac{P_{i,j}}{e_{1}}}\otimes e_{i}\\
\partial_{3}^{h}(1 \otimes e_{i}\wedge e_{j}\wedge e_{k})=
\partial_{3}(1 \otimes e_{i}\wedge e_{j}\wedge e_{k})-
{\displaystyle \mathop \frac{P_{i,j}}{e_{1}}}\otimes e_{1}\wedge e_{k}
+{\displaystyle \mathop \frac{P_{i,k}}{e_{1}}}\otimes e_{1}\wedge
e_{j}-
{\displaystyle \mathop \frac{P_{j,k}}{e_{1}}}\otimes e_{1}\wedge
e_{i}\\
\partial_{4}^{h}(1 \otimes e_{1}\wedge e_{i}\wedge e_{j}\wedge e_{k})=
\partial_{4}(1 \otimes e_{1}\wedge e_{i}\wedge e_{j}\wedge e_{k})\\
\partial_{4}^{h}(1 \otimes e_{2}\wedge e_{3}\wedge e_{4}\wedge e_{5})=
\partial_{4}(1 \otimes e_{2}\wedge e_{3}\wedge e_{4}\wedge e_{5})
+{\displaystyle \mathop \frac{P_{3,5}}{e_{1}}}
\otimes e_{1}\wedge e_{2}\wedge e_{4} \\
\;\;\;\;\;\;\;\;\;\;\;\;\;\;-{\displaystyle \mathop \frac{P_{3,4}}{e_{1}}}
\otimes e_{1}\wedge e_{2}\wedge e_{5}
-{\displaystyle \mathop \frac{P_{4,5}}{e_{1}}}
\otimes e_{1}\wedge e_{2}\wedge e_{3}
-{\displaystyle \mathop \frac{P_{2,5}}{e_{1}}}
\otimes e_{1}\wedge e_{3}\wedge e_{4}\\
\partial_{5}^{h}(1 \otimes 
e_{1}\wedge e_{2}\wedge e_{3}\wedge e_{4}\wedge e_{5})=
\partial_{5}(1 \otimes 
e_{1}\wedge e_{2}\wedge e_{3}\wedge e_{4}\wedge e_{5}).\\
\end{array}$$
The character defined by the right multiplication of 
$U_{h}({\goth a})$ on 
$Ext^{5}_{U_{h}({\goth a})}\left ( k[[h]], U_{h}({\goth a})\right )$
is zero. 
\end{proposition}
{\it Proof of the proposition : }
The resolution of $k[[h]]$ constructed in the proposition is obtained
by applying the proof of theorem ~\ref{koszul resolution}. Moreover, one has

$$\begin{array}{l}
^{t}\partial_{5}(1\otimes e_{1}^{*}\wedge e_{2}^{*}\wedge
e_{3}^{*}\wedge e_{4}^{*})=e_{5}\otimes 
e_{1}^{*}\wedge e_{2}^{*}\wedge
e_{3}^{*}\wedge e_{4}^{*}\wedge e_{5}^{*}\\
 ^{t}\partial_{5}(1\otimes e_{1}^{*}\wedge e_{3}^{*}\wedge
e_{4}^{*}\wedge e_{5}^{*})=-e_{2}\otimes 
e_{1}^{*}\wedge e_{2}^{*}\wedge
e_{3}^{*}\wedge e_{4}^{*}\wedge e_{5}^{*}\\
^{t}\partial_{5}(1\otimes e_{1}^{*}\wedge e_{2}^{*}\wedge
e_{4}^{*}\wedge e_{5}^{*})=e_{3}\otimes e_{1}^{*}\wedge e_{2}^{*}\wedge
e_{3}^{*}\wedge e_{4}^{*}\wedge e_{5}^{*}\\
^{t}\partial_{5}(1\otimes e_{1}^{*}\wedge e_{2}^{*}\wedge
e_{3}^{*}\wedge e_{5}^{*})=-e_{4}\otimes e_{2}^{*}\wedge e_{3}^{*}\wedge
e_{3}^{*}\wedge e_{4}^{*}\wedge e_{5}^{*}\\
^{t}\partial_{5}(1\otimes e_{2}^{*}\wedge e_{3}^{*}\wedge
e_{4}^{*}\wedge e_{5}^{*})=e_{1}\otimes e_{1}^{*}\wedge e_{2}^{*}\wedge
e_{3}^{*}\wedge e_{4}^{*}\wedge e_{5}^{*}.\\ 
\end{array}$$
These equalities show that the  character defined by the right 
multiplication of $U_{h}({\goth a})$ on 
$Ext^{5}_{U_{h}({\goth a})}\left ( k[[h]], U_{h}({\goth a})\right )$
is zero. 

\section{Applications}

\subsection{Poincar{\'e} duality}

Let $M$ be an $A_{h}^{op}$-module and $N$ an 
$A_{h}$-module. The right exact functor 
$M {\displaystyle  \mathop \otimes_{A_{h}}}-$ has
a left derived functor. We set   
$Tor^{i}_{A_{h}} \left (M, N \right )=
L^{i}\left ( M  {\displaystyle  \mathop \otimes_{A_{h}}}-\right ) (N)$. \\

\begin{theorem}~\label{Poincare}
Let $A_{h}$ be a deformation algebra of $A_{0}$ satisfying the
hypothesis of theorem ~\ref{qtrad}.  Assume moreover that the $A_{h}$-module $K$ is   of finite projective dimension. 
Let $M$ be an $A_{h}$-module. 
One has an isomorphism of $K$-modules 
$$Ext^{i}_{A_{h}}\left ( K, M \right ) \simeq 
Tor^{A_{h}}_{d_{A_{h}}-i}
\left (\Omega_{A_{h}}, M \right ) .$$
\end{theorem} 

{\bf Remark : } Theorem ~\ref{Poincare} generalizes classical Poincar{\'e}
duality  ([Kn]). \\

{\it Proof of the theorem}

As the $A_{h}$-module $K$ admits a finite length resolution by finitely
generated  projective $A_{h}$-modules, 
$P^{\bullet} \to K$, the
canonical arrow 
$$RHom_{A_{h}}\left ( K, A_{h}\right )
{\displaystyle \mathop \otimes_{A_{h}}^{L}}M \to 
RHom _{A_{h}}(K,M )$$ is an isomorphism in $D(Mod A_{h})$. Indeed the canonical
arrow 
$$Hom_{A_{h}}\left ( P^{\bullet}, A_{h}\right )
{\displaystyle \mathop \otimes_{A_{h}}}M \to 
Hom _{A_{h}}(P^{\bullet},M )$$
is an isomorphism. 

\subsection{Duality property for induced representations of quantum
  groups}$\;$\\

 {\it From now on, we assume that $A_{h}$ is a topological Hopf algebra.}\\

 In this section, we keep the notation of theorem \ref{two dualities}. 
 Let $V$ be a left $A_{h}$-module, then, by transposition, 
 $V^{*}=Hom_{K}(V,K)$ is naturally endowed with a  right 
$A_{h}$-module structure. Using the antipode, we can also see $V^{*}$ as a
left module structure. Thus, one has : 
$$\forall u \in A_{h}\, \forall f \in V^{*}, \;\; u \cdot f =f \cdot S(u).$$
We endow $ \Omega_{A_{h}}\otimes V^{*}$ with the following
right  $A_{h}$-module structure :
$$\begin{array}{l}
\forall u \in A_{h}\, \forall f \in V^{*}, \;
\forall \omega \in \Omega_{A_{h}},\\
(\omega \otimes f)\cdot u =
 \lim\limits_{n\to +\infty}
{\displaystyle \mathop \sum_{j}}
\theta_{A_{h}} (u'_{j,n})\omega  \otimes 
f \cdot S_{h}^{2}(u''_{j,n}) 
\end{array}$$
where $\Delta (u)=\lim\limits_{n\to +\infty}
{\displaystyle \mathop \sum_{j}u'_{j,n}\otimes u''_{j,n}}$.

Let  $A_{h}$ be a topological Hopf deformation of $A_{0}$ and $B_{h}$ be a
topological Hopf deformation of $B_{0}$. We assume moreover that there
exists a morphism of Hopf algebras from $B_{h}$ to $A_{h}$ and 
that $A_{h}$ is a flat $B_{h}^{op}$-module (by proposition ~\ref{flatness}
this is verified if 
the induced $B_{0}$-module structure on $A_{0}$ is flat). If $V$ is an
 $A_{h}$-module, we can define the induced 
representation from $V$ as follows :
$$Ind_{B_{h}}^{A_{h}}\left ( V\right )=
A_{h}
{\displaystyle \mathop \otimes_{B_{h}}}V$$
on which $A_{h}$ acts by left multiplication.\\

%{\bf Remark :}

%Considered as a right $B_{h}$-module using the antipode, $V$ will be denoted
%$V^{r}$. We will write 
%$\left ( V^{r} {\displaystyle \mathop \otimes_{B_{h}}}A_{h}\right )^{l}$ 
%for the  $V^{r} {\displaystyle \mathop \otimes_{B_{h}}}A_{h}$
%considered as a left $A_{h}$-module. 
%The map 
%$\Psi : A_{h}{\displaystyle \mathop \otimes_{B_{h}}}V \to  
%\left ( V^{r} {\displaystyle \mathop \otimes_{B_{h}}}A_{h}\right
%)^{l}$ defined by 
%$$\forall a \in A_{b}, \forall v \in V, \;\; \Psi (a \otimes v)=v
%\otimes S(a)$$
%is an  isomorphism of left $A_{h}$-modules. 
%In the sequel, we will omit the upperscripts
%$l$ and $r$. \\

\begin{proposition}~\label{Ext-dual}
Let  $A_{h}$ be a topological Hopf deformation of $A_{0}$ and $B_{h}$ be a
topological deformation of $B_{0}$. We assume that there
exists a morphism of Hopf algebras from $B_{h}$ to $A_{h}$ such that
$A_{h}$ is a flat  $B_{h}^{op}$-module. We also
assume that $B_{h}$ satisfies the hypothesis of theorem ~\ref{qtrad}. 
Let $V$ be an 
$B_{h}$-module which is a  free finite dimensional K-module. 
Then 
$D_{B_{h}}\left ( Ind^{B_{h}}_{A_{h}}(V)
  \right )$ is isomorphic to 
$ \left ( \Omega_{B_{h}} \otimes V^{*}\right )_
{\displaystyle \mathop \otimes _{B_{h}}}A_{h}[-d_{B_{h}}]$ in 
$D\left (Mod B_{h}^{op} \right )$.\\
\end{proposition}

\begin{corollary}~\label{corollaire Ext-dual}
Let $A_{h}$ be a topological Hopf deformation of $A_{0}$ and $B_{h}$ be a
topological deformation of $B_{0}$. We assume that there exists a
morphism of Hopf algebras from $B_{h}$ to $A_{h}$ such that $A_{h}$ is
a flat $B_{h}^{op}$-module. We also assume that $B_{h}$ satisfies the
condition of the theorem ~\ref{qtrad}. Let $V$ be a $B_{h}$-module which is
a free finite dimensional $K$-module. Then 

a) $Ext^{i}_{A_{h}} \left ( 
A_{h}{\displaystyle \mathop \otimes_{B_{h}}}V, A_{h} \right )$ is
reduced to $0$ if $i$ is different from $d_{B_{h}}$.

b) The right $A_{h}$-module 
 $Ext^{d_{B_{h}}}_{A_{h}} \left ( 
A_{h}{\displaystyle \mathop \otimes_{B_{h}}}V, A_{h} \right )$ is
isomorphic to 
$\left (  \Omega_{B_{h}}\otimes V^{*}\right ){\displaystyle \mathop \otimes_{B_{h}}}A_{h}$.
\end{corollary}

{\bf Remarks :}

Proposition ~\ref{Ext-dual} is already known in
the case where $\mathfrak g$ is a Lie algebra, ${\mathfrak h}$  
is a Lie subalgebras of ${\mathfrak g}$, 
$A$ and  $B$ are the corresponding enveloping algebras. 
 
In this case one has $d_{B_{h}}=dim{\mathfrak h}$ and 
$d_{C_{h}}=dim{\mathfrak k}$. More
 precisely : It was proved by Brown
and Levasseur ([B-L] p. 410) and [Ke] in the case where $\mathfrak g$ is a
finite dimensional semi-simple Lie algebra and 
$Ind_{U({\mathfrak h})}^{U({\mathfrak g})}(V)$ is a Verma-module. 
Proposition ~\ref{duality property}
is proved in full generality for Lie superalgebras in [C1].\\

Here are some examples of  situations where we can apply the
proposition ~\ref{Ext-dual}:\\

{\it Example 1 :}

Let $k$ be a field of characteristic $0$. We set $K=k[[h]]$. Etingof
and Kazhdan  have constructed a functor $Q$ from the category $LB(k)$ of Lie
bialgebras over $k$ to the category $HA(K)$ of topological Hopf
algebras over $K$. If $({\mathfrak g}, \delta)$ is a Lie bialgebra, its
image by $Q$ will be denoted $U_{h}({\mathfrak g})$. 

Let $\mathfrak g$ be a Lie bialgebra 
Let $\mathfrak h$ be a Lie sub-bialgebra
of ${\mathfrak g}$. The functoriality of the quantization
implies  the existence of an embedding  of Hopf algebras
from $U_{h}({\mathfrak h})$ to $U_{h}({\mathfrak g})$ which satisfies all our
hypothesis.\\

{\it Example 2 :} If $\mathfrak g$ is a Lie bialgebra, we will denote
by ${\mathcal F}({\mathfrak g})$ the formal group attached to it and 
${\mathcal F}_{h}({\mathfrak g})$ its Etingof Kazhdan quantization. 
Let ${\mathfrak g}$ and $\mathfrak h$ be two Lie algebras and assume
that there exists a surjective morphism of Lie bialgebras 
from ${\mathfrak g}$ to ${\mathfrak h}$. Then 
${\mathcal F}_{h}({\mathfrak g})$ is a flat 
${\mathcal F}_{h}({\mathfrak h})$-module and 
$A_{h}={\mathcal F}_{h}({\mathfrak g})$
and $B_{h}={\mathcal F}_{h}({\mathfrak h})$ 
satisfies the hypothesis of the theorem.\\  

{\it Example 3 :}

If $G$ is an affine algebraic Poisson  group, we will denote by 
${\mathcal F}(G)$ the algebra of regular functions on $G$ and 
${\mathcal F}_{h}(G)$ its Etingof Kazhdan quantization. 
Let $G$ and $H$ be  affine algebraic Poisson  groups. Assume  that 
there is a Poisson  group map $G \to H$ 
such that ${\mathcal F}(G)$ is a flat ${\mathcal F}(H)^{op}$-module. 
By functoriality of Etingof Kazhdan quantization, $A_{h}={\mathcal F}_{h}(G)$
and $B_{h}={\mathcal F}_{h}(H)$ satisfies the hypothesis of the
theorem.  \\

{\it Proof of the proposition ~\ref{Ext-dual} :}

We proceed as in [C1]. Let $L^{\bullet}\to V$ be a resolution of $V$ by finite free 
$B_{h}$-modules. As $A_{h}$ is a flat $B_{h}^{op}$-module, 
$A_{h}{\displaystyle \mathop \otimes_{B_{h}}}L^{\bullet} \to 
A_{h}{\displaystyle \mathop \otimes_{B_{h}}}V$ is a resolution of the $A_{h}$-module 
$A_{h}{\displaystyle \mathop \otimes_{B_{h}}}V$ by finite free $A_{h}$-modules. 

We have the following sequence of isomorphisms in 
$D\left ( Mod A_{h} \right )$
$$\begin{array}{rcl}
RHom_{A_{h}}\left ( 
A_{h}{\displaystyle \mathop \otimes_{B_{h}}}V,
A_{h} \right ) & \simeq & 
Hom_{A_{h}}\left ( 
A_{h}{\displaystyle \mathop \otimes_{B_{h}}}L^{\bullet},
A_{h} \right )\\
& \simeq & 
Hom_{B_{h}}(L^{\bullet}, B_{h}){\displaystyle \mathop \otimes_{B_{h}}}A_{h}\\
&\simeq & 
\left ( \Omega_{B_{h}}\otimes V^{*} \right ){\displaystyle \mathop \otimes_{B_{h}}}A_{h}
[-d_{B_{h}}]. \Box
\end{array}$$ \\

We now extend to Hopf algebras another duality property for induced
representations of Lie algebras ([C1]).\\  

\begin{proposition}~\label{duality property}
Let  $A_{h}$ be a Hopf deformation of $A_{0}$,  $B_{h}$ be a
Hopf  deformation of $B_{0}$ and  
$C_{h}$ be a Hopf  deformation of $C_{0}$. We assume that there
exists a morphism of Hopf algebras from $B_{h}$ to $A_{h}$ 
and a morphism of Hopf algebras from $C_{h}$ to $A_{h}$ 
such that $A_{h}$
is a flat  $B_{h}^{op}$-module and a flat $C_{h}^{op}$-module. We also
assume that $B_{h}$ and $C_{h}$ satisfies the hypothesis of 
theorem ~\ref{qtrad}. 
Let $V$ (respectively $W$) be an 
$B_{h}$-module (respectively $C_{h}$-module) which is a  free finite 
dimensional K-module. Then, for all integer $n$, one has an
isomorphism 

$$\begin{array}{l}
Ext^{n +d_{B_{h}}}_{A_{h}}\left ( 
A_{h}{\displaystyle \mathop \otimes_{B_{h}}} V, 
A_{h}{\displaystyle \mathop \otimes_{C_{h}}} W
\right )\\
\simeq 
Ext^{n +d_{C_{h}}}_{A_{h}^{op}}\left ( 
\left (  \Omega_{C_{h}}\otimes W^{*}\right ){\displaystyle \mathop \otimes_{C_{h}}}A_{h} , 
\left (  \Omega_{B_{h}}\otimes V^{*}\right ){\displaystyle \mathop \otimes_{C_{h}}} A_{h}\right )
\end{array} $$
\end{proposition}

{\bf Remarks : } 

Proposition ~\ref{duality property} is already known in
the case where $\mathfrak g$ is a Lie algebra, ${\mathfrak h}$ and 
${\mathfrak k}$ are Lie subalgebras of ${\mathfrak g}$, 
$A$, $B$ and  $C$ are the corresponding enveloping algebras. 
In this case one has $d_{B_{h}}=dim{\mathfrak h}$ and 
$d_{C_{h}}=dim{\mathfrak k}$. More
 precisely : 

 Generalizing a result of G. Zuckerman ([B-C]), A. Gyoja ([G]) proved a
 part of this theorem (namely the case where 
${\mathfrak h}={\mathfrak g}$ and $n=dim{\mathfrak h}=dim {\mathfrak
  k}$) under the assumptions that ${\mathfrak g}$ is split semi-simple
and ${\mathfrak h}$ is a parabolic subalgebra of ${\mathfrak g}$. 
D.H Collingwood and B. Shelton ([C-S]) also proved a duality of this
type (still under the semi-simple hypothesis) 
but in a slighly different context.  

M. Duflo [Du2] proved proposition ~\ref{duality property} for a
${\mathfrak g}$ general Lie algebra, ${\mathfrak h}={\mathfrak k}$,
$V=W^{*}$ being one dimensional representations.

Proposition ~\ref{duality property} is proved in full generality in the
context of Lie superalgebras in [C1].\\

{\it Proof of the proposition  ~\ref{duality property}:}

We will proceed as in [C2]. 
%As $D$ is a contravariant functor from $D^{b}_{f}(A_{h})$ to itself and 
As $D_{A_{h}^{op}}\circ D_{A_{h}}\left ( 
A_{h}{\displaystyle \mathop \otimes_{B_{h}}} V
\right )= 
A_{h}{\displaystyle \mathop \otimes_{B_{h}}} V$,
we have the following isomorphism 

$$\begin{array}{l}
Hom_{D(A_{h})}\left ( 
A_{h}{\displaystyle \mathop \otimes_{B_{h}}} V, 
A_{h}{\displaystyle \mathop \otimes_{B_{h}}} W
\right )\\
\simeq
Hom_{D(A_{h}^{op})}\left [  
D_{A_{h}}\left ( 
A_{h}{\displaystyle \mathop \otimes_{C_{h}}} W
\right ), 
D_{A_{h}}\left ( 
A_{h}{\displaystyle \mathop \otimes_{B_{h}}} V 
\right ) \right ]
\end{array}$$

the corollary follows now from proposition ~\ref{Ext-dual}. 

\subsection{Hochschild cohomology} In this subsection, $A_{h}$ is a topological Hopf algebra. 
We set $A_{h}^{e}=A_{h}{\displaystyle \mathop \otimes _{k[[h]]}A_{h}^{op}}$ and 
$\widehat{A_{h}^{e}}=
A_{h}\widehat{\displaystyle \mathop \otimes _{k[[h]]}}A_{h}^{op}$. If $M$ is an  
$\widehat{A_{h}^{e}}$-module, we set 
$$\begin{array}{l}
HH^{i}_{A_{h}}(M)= Ext^{i}_{\widehat{A_{h}^{e}}}(A_{h}, M)\\
HH_{i}^{A_{h}}(M)=Tor^{\widehat{A_{h}^{e}}}_{i}(A_{h}, M)
\end{array}$$

\begin{proposition}\label{Hochschild}
Assume that $A_{h}$ satisfies the condition of the theorem \ref{qtrad}. Assume moreover that 
$A_{0}\otimes A_{0}^{op}$ is noetherian. Consider 
$A_{h}\widehat{\displaystyle \mathop \otimes _{k[[h]]}}A_{h}$ with the following 
$\widehat{A_{h}^{e}}$-module structure :
$$\forall (\alpha , \beta , x ,y ) \in A_{h}, \;\;\;\alpha \cdot (x \otimes y) \cdot \beta = 
\alpha x \otimes y \beta .$$

a) $HH^{i}_{A_{h}}(A_{h}
\widehat{\displaystyle \mathop \otimes _{k[[h]]}}A_{h})$ is zero if $i\neq d_{A_{h}}$.

b)  The $\widehat{A_{h}^{e}}$-module  $HH^{d_{A_{h}}}_{A_{h}}
(A_{h}\widehat{\displaystyle \mathop \otimes _{k[[h]]}}A_{h})$ is isomorphic to 
$\Omega_{A_{h}} \otimes A_{h}$ with the following $\widehat{A_{h}^{e}}$-module structure : 
$$\forall (\alpha , \beta , x ) \in A_{h}, \;\;\;\ \\
\alpha \cdot (\omega \otimes x) \cdot \beta = 
\omega \theta_{A_{h}}(\beta '_{i})  \otimes S(\beta_{i}'')
 x S^{-1}( \alpha) $$
 where $\alpha={\displaystyle \sum_{i}}\alpha_{i}'\otimes \alpha''_{i}$ (to be taken in the topological sense)
\end{proposition}
{\it Proof of the theorem :}

The proof is analogous to that of [C2] (theorem 3.3.2). 

Using the antipode $S_{h}$ of $A_{h}$, we have the following isomorphism  in 
$D \left ( Mod \widehat {A_{h}^{e}}\right ) $, 
$$RHom_{\widehat{A_{h}^{e}}}\left ( A_{h}, 
A_{h}\widehat {\displaystyle \mathop \otimes} A_{h}\right ) \simeq 
RHom_{A_{h} \widehat {\displaystyle \mathop \otimes} A_{h}}\left ( (A_{h})^{\#}, 
(A_{h} \widehat {\displaystyle \mathop \otimes} A_{h})^{\#}\right ) .$$
where the structures on $(A_{h})^{\#}$ and 
$(A_{h} \widehat {\displaystyle \mathop \otimes} A_{h})^{\#}$ are given by : 
$$\begin{array}{l}
\forall (\alpha , \beta , u,v) \in A_{h}\\
(\alpha \otimes \beta ) \cdot u = \alpha u S_{h}(\beta )\\
(\alpha \otimes \beta )\cdot (u \otimes v)= \alpha u \otimes v S_{h}(\beta)\\
(u \otimes v ) \cdot \alpha \otimes \beta = u \alpha \otimes S_{h}(\beta )v. 
\end{array}$$
Using the version of  lemma \ref{bimodule isomorphism} for right modules (see [C2] lemma 1;1),  one sees that $ (A_{h})^{\#}$ is isomophic to 
$(A_{h} \widehat {\displaystyle \mathop \otimes} A_{h})
{\displaystyle \mathop \otimes_{A_{h}}}K$ 
as an $A_{h}\widehat {\displaystyle \mathop \otimes} A_{h}$-module. 
we get 
$$\begin{array}{rcl}
RHom_{\widehat{A_{h}^{e}}}\left ( A_{h}, 
A_{h}\widehat {\displaystyle \mathop \otimes} A_{h}\right )& \simeq &
RHom_{A_{h} \widehat {\displaystyle \mathop \otimes} A_{h}}\left ( 
A_{h} \widehat {\displaystyle \mathop \otimes} A_{h}
{\displaystyle \mathop \otimes_{A_{h}}}K, 
(A_{h} \widehat {\displaystyle \mathop \otimes} A_{h})^{\#}\right ) \\
& \simeq & 
RHom_{ A_{h}}\left ( K, 
(A_{h} \widehat {\displaystyle \mathop \otimes} A_{h})^{\#}\right ) \\
& \simeq & 
RHom_{ A_{h}}\left ( K,  A_{h}\right )
  {\displaystyle \mathop \otimes}_{A_{h}} 
 (A_{h}  \widehat{\displaystyle \mathop \otimes} A_{h})^{\#}\\
 & \simeq & 
\Omega_{h}
  {\displaystyle \mathop \otimes}_{A_{h}} 
 (A_{h} \widehat {\displaystyle \mathop \otimes} A_{h})^{\#}
 \end{array} $$
 The isomorphism $id \otimes S_{h}^{-1}$ transforms 
 $(A_{h} \widehat {\displaystyle \mathop \otimes} A_{h})^{\#}$ into the natural 
$(A_{h} \widehat {\displaystyle \mathop \otimes} A_{h})
\otimes (A_{h} \widehat {\displaystyle \mathop \otimes} A_{h})^{op}$-module 
$(A_{h} \widehat {\displaystyle \mathop \otimes} A_{h})$-module 
$(A_{h} \widehat {\displaystyle \mathop \otimes} A_{h})^{nat}$:  
$$\begin{array}{l}
\forall (\alpha , \beta , u,v) \in A_{h}\\
(\alpha \otimes \beta )\cdot (u \otimes v)= \alpha u \otimes \beta  v \\
(u \otimes v ) \cdot \alpha \otimes \beta = u \alpha \otimes v \beta. 
\end{array}$$
Then, 
using the  lemma \ref{bimodule isomorphism},  one sees that 
$\Omega_{h}
  {\displaystyle \mathop \otimes}_{A_{h}} 
 (A_{h} \widehat {\displaystyle \mathop \otimes} A_{h})^{nat}$ is isomorphic to 
 $\Omega_{h}\otimes A_{h}$ endowed with the following 
$ (A_{h} \widehat {\displaystyle \mathop \otimes} A_{h})^{op}$-module structure :
$$
\forall (\alpha , \beta ) \in A_{h}, \;\; \\
(u \otimes v ) \cdot \alpha \otimes \beta = 
{\displaystyle \sum_{i}}u \theta_{A_{h}}(\alpha_{i}') \otimes 
S(\alpha_{i}'')v \beta .$$
This finishes the proof of the proposition. $\Box$\\

We are in the case where $Ext^{i}_{\widehat{A_{h}^{op}}}(A_{h}, \widehat{A_{h}^{e}})$ is $0$ unless when 
$i=d_{A_{h}}$, so we have a duality between Hochschild homology and Hochschild cohomology ([VdB]). 

\begin{corollary}
Let $A_{h}$ be a $k$-algebra satisfying the hypothesis of theorem \ref{qtrad}. Assume moreover that 
$A_{0}^{e}=A_{0}\otimes A_{0}^{op}$ is noetherian and that the $\widehat{A_{h}^{e}}$-module 
$A_{h}$ is finite projective dimension. Let $M$ be an $\widehat{A_{h}^{e}}$-module. 
One has 
$$HH^{i}(M)\simeq HH_{d_{A_{h}}-i}(U {\displaystyle \mathop \otimes_{A_{h}}} M).$$
\end{corollary}

{\it Proof of the corollary :} The proof of the corollary is similar to that of [vdB].

{\it First case : $M$ is a finite type $\widehat{A_{h}^{e}}$-module.} Let 
$P^{\bullet}\to A_{h}\to 0$ be a finite length finite type projective resolution of the  
$\widehat{A_{h}^{e}}$-module $A_{h}$ and let $Q^{\bullet}\to M \to 0$ be a finite type projective resolution of the $\widehat{A_{h}^{e}}$-module $M$. 
As $Q^{i}$ and $U{\displaystyle  \mathop \otimes _{A_{h}}}Q^{i}$ are complete, 
one has the following sequence of isomorphisms :
$$\begin{array}{rcl}
HH^{i}_{\widehat{A_{h}^{e}}}(M)
&\simeq &H^{i}\left ( Hom_{\widehat{A_{h}^{e}}}(P^{\bullet}, M)\right )
\simeq H^{i}\left ( Hom_{\widehat{A_{h}^{e}}}(P^{\bullet}, \widehat{A_{h}^{e}})
{\displaystyle  \mathop \otimes _{\widehat{A_{h}^{e}}}}M\right )\\
&\simeq &H^{i}\left ( U[-d]
{\displaystyle  \mathop \otimes^{L} _{\widehat{A_{h}^{e}}}}M\right )
\simeq 
H^{i-d_{A_{h}}}\left ( U{\displaystyle  \mathop \otimes _{\widehat{A_{h}^{e}}}}Q^{\bullet}\right )
\simeq 
H^{i-d_{A_{h}}}\left ( (A_{h}{\displaystyle  \mathop \otimes _{A_{h}}}
U){\displaystyle  \mathop \otimes _{\widehat{A_{h}^{e}}}}Q^{\bullet})\right )\\
&\simeq &
H^{i-d_{A_{h}}}\left ( A_{h}{\displaystyle  \mathop \otimes _{\widehat{A_{h}^{e}}}}
(U{\displaystyle  \mathop \otimes _{A_{h}}}Q^{\bullet})\right )
\simeq HH_{d_{A_{h}}-i}(U 
{\displaystyle  \mathop \otimes _{A_{h}}} M).
\end{array}$$

{\it General case :} We no longer assume that $M$ is a finite type $\widehat{A_{h}^{e}}$-module. 
We have $M=\lim\limits_{\rightarrow}M'$ 
where $M'$ runs over all finitely generated submodules of $M$. 
$$\begin{array}{l}
\Ext^{i}_{\widehat{A_{h}^{e}}}(A_{h},M)= 
\Ext^{i}_{\widehat{A_{h}^{e}}}(A_{h},\lim\limits_{\rightarrow}M')\simeq 
\lim\limits_{\rightarrow}\Ext^{i}_{\widehat{A_{h}^{e}}}(A_{h},M')\simeq
\lim\limits_{\rightarrow}  
Tor_{d_{A_{h}}-i}^ {\widehat{A_{h}^{e}}}
(A_{h}, U {\displaystyle  \mathop \otimes _{\widehat{A_{h}}}} M')\\ \simeq 
Tor_{d_{A_{h}}-i}^{\widehat{A_{h}^{e}}} (A_{h}, \lim\limits_{\rightarrow}
U {\displaystyle  \mathop \otimes _{\widehat{A_{h}}}} M') 
\simeq Tor_{d_{A_{h}}-i}^ {\widehat{A_{h}^{e}}}(A_{h}, 
U {\displaystyle  \mathop \otimes _{\widehat{A_{h}}}} M) 
\end{array}$$
where  we used the fact that
the functor $\lim\limits_{\rightarrow}$ is exact because the set of
finitely generated submodules of $M$ is a directed set 
([Ro] proposition 5.33)

\section{Bibliography}

[A-K] A.Altman and S. Kleiman : {\it Introduction to Grothendieck
  duality theory}, Lecture Notes in Mathematics {\bf 146}, Springer-Verlag. 

[B-C] B.D Boe and D.H Collingwood : {\it A comparison theorem for the
  structures of induced representations }, J. algebra {\bf 94}, 1985, 
p. 511-545.

[B-L] K. A. Brown - T. Levasseur : {\it Cohomology of bimodules over
  enveloping algebras}, Mathematische Zeitschrift {\bf 189}, 1985,
511-545. 

[C-P] V. Chari - A. Pressley : {\it A guide to quantum groups},
Cambridge University Press (1994).

[C1] S. Chemla : {\it Poincare duality for k-A Lie superalgebras}, 
Bulletin de la Soci{\'e}t{\'e} Math{\'e}matique de
France {\bf 122}, 371-397. 

[C2] S. Chemla : {\it Rigid dualizing complex for quantum algebras
  and algebras of generalized differential operators}, Journal of
algebra {\bf 276} (2004) p.80-102.

[Chl] Chloup-Arnoult V. :{\it Linearization of some Poisson-Lie
  tensor}, Journal of geometry and physics {\bf 24} (1997) 46-52.

[C-S] D.H Collingwood- B. Shelton :{\it A duality theorem for
  extensions of induced highest weight modules}, Pacific J. Math {\bf
  146} 2, 1990, 227-237. 

[Dr] V. G. Drindfeld : {\it Quantum groups}, Proc. Intern. Congress of
Math. (Berkeley 1986)(1987), 798-820. 

[Du1] M. Duflo : 
{\it Sur les id{\'e}aux induits dans les alg{\`e}bres enveloppantes}
Invent. Math. {\bf 67}, 385-393.

[Du2] M. Duflo : {\it Open problems in representation theory of Lie
  groups}, Proceedings of the eighteenth international symposium,
division of mathematics, the Tanigushi foundation. 

[E-K 1] P. Etingof - Kazhdan :{\it Quantization of Lie bialgebras I},
Selecta mathematica{\bf 2} (1996) p1-41.

[E-K 2] P. Etingof - Kazhdan :{\it Quantization of Lie bialgebras II},
Selecta mathematica. 

[E-K 3] P. Etingof - Kazhdan :{\it Quantization of Poisson algebraic
  groups and Poisson homogeneous spaces}, ArXiv:q-alg/9510020v2.

[E-S] P. Etingof - O. Schiffman : {\it Lectures on quantum groups}, 
International Press of Boston Inc (2001). 

[Ga] F. Gavarini : {\it The quantum duality principle},
Ann. Inst. Fourier, {\bf 52}, 3 (2002), 809-834.

[Gy] A. Gyoja : 
{\it A duality theorem for homomorphisms between generalized Verma modules} 
Preprint Kyoto University.

[H-S] P.J Hilton - U. Stammbach {\it A course in homological algebra}, 
Graduate text in mathematics, Springer Verlag. 

[K-S] M. Kashiwara - P. Schapira : {\it Deformation quantization
  modules I} arXiv : math QA/0802.1245.

[Ke] G.R. Kempf : {\it The Ext-dual of a Verma module is a Verma
  module}, Journal of pure and applied algebra, {\bf 75}, 1991,
p. 47-49. 

[Kn] A. Knapp : {\it Lie groups, Lie algebras and cohomology},
Princeton university Press, 1988.

[L] J.L. Loday : {\it Cyclic homology}, Springer Verlag, Berlin
(1998). 

[Ro] J.J Rotman : {\it An introduction to homological algebra}, second
edition, Springer.

[Schn] J-P Schneiders : {\it An introduction to D-modules}, Bulletin
Soc. Royale Science Liege {\bf 63} (3-4)(1994) p. 223-295.

[Schw] L. Schwartz : {\it Topologie g{\'e}n{\'e}rale et analyse fonctionnelle},
Hermann. 

[VdB] M. Van den Bergh : {\it A relation between Hochschild homology and cohomology for Gorenstein rings}, Proc. Amer. Math. Soc. {\bf 126} (1998), 1345-1348 and erratum in Proc. Amer. Math. Soc.  
{\bf 130} (2002), 2809-2810. 
%[Ye] A. Yekutieli : {\it The rigid dualizing complex of a universal
%  enveloping algebra}, Journal of pure and applied algebra{\bf 150}
%(2000) p.85-93.

\end{document}